\documentclass[11pt,twoside]{article}
\usepackage[T1]{fontenc}
\usepackage{amsmath,amsfonts,amssymb,amsthm,fullpage,bbm}
\usepackage{fancyhdr,graphicx,color}
\usepackage{epsfig}
\usepackage{yhmath}

\newtheorem{thm}{Theorem}

\newtheorem{lem}{Lemma}

\theoremstyle{definition}

\theoremstyle{remark}
\newtheorem{rem}{Remark}

\DeclareMathOperator{\cyl}{cyl}
\DeclareMathOperator{\card}{card}
\DeclareMathOperator{\di}{div}
\DeclareMathOperator{\disc}{disc}
\DeclareMathOperator{\hyp}{hyp}
\DeclareMathOperator{\diam}{diam}

\DeclareMathOperator{\flow}{flow}

\newcommand{\eps}{\varepsilon}

\def\PP{\mathbb{P}}
\def\RR{\mathbb{R}}
\def\EE{\mathbb{E}}

\def\ZZ{\mathbb{Z}}
\def\E{\mathcal{E}}
\def\D{\mathcal{D}}
\def\H{\mathcal{H}}

\def\I{\mathcal{I}}

\def\G{\Gamma}
\def\O{\Omega}
\def\p{\partial}
\def\L{\mathcal{L}}
\def\P{\mathcal{P}}
\def\B{\mathcal{B}}
\def\C{\mathcal{C}}

\def\R{\mathcal{R}}
\def\V{\mathcal{V}}
\def\U{\mathcal{U}}

\def\ind{{\mathbbm{1}}_}

\newcommand{\bro}{\smash{{B}^{\!\!\!\!\raise5pt\hbox{\scriptsize o}}}}
\newcommand{\dro}{\smash{{D}^{\!\!\!\!\raise5pt\hbox{\scriptsize o}}}}
\newcommand{\lro}{\smash{{L}^{\!\!\!\!\raise5pt\hbox{\scriptsize o}}}}
\newcommand{\qro}{\smash{{Q}^{\!\!\!\!\raise5pt\hbox{\scriptsize o}}}}
\newcommand{\uro}{\smash{{U}^{\!\!\!\!\raise5pt\hbox{\scriptsize o}}}}
\newcommand{\ero}{\smash{{E}^{\!\!\!\!\raise5pt\hbox{\scriptsize o}}}}
\newcommand{\cero}{\smash{{\calE}^{\!\!\!\!\raise5pt\hbox{\scriptsize o}}}}
\newcommand{\cgcirc}{\smash{{\calG}^{\!\!\!\raise5pt\hbox{\scriptsize o}}}}
\newcommand{\cdero}{\smash{{\mathbb{E}}^{\!\!\!\!\raise5pt\hbox{\scriptsize o}}}}
\newcommand{\bgcirc}{\smash{{\mathbb{G}}^{\!\!\!\!\raise5pt\hbox{\scriptsize o}}}}
\newcommand{\begcirc}{{{{\mathbb E}\cap\GGG}^{
\!\!\!\!\!\!\!\!\!\raise3pt\hbox{\scriptsize o}}}}

\newcommand{\aro}{\smash{{A}^{\!\!\!\raise5pt\hbox{\scriptsize o}}}}
\newcommand{\thro}{\smash{{\Th}^{\!\!\!\!\!\raise5.2pt\hbox{\scriptsize
        o}}}}

\title{\huge Lower large deviations for the maximal flow through a domain
  of $\RR^d$ in first passage percolation} 
\author{}
\date{}

\begin{document}
\maketitle

\thispagestyle{empty}

\begin{center}
\vskip-1cm {\Large Rapha\"el Cerf}\\
{\it Universit\'e Paris Sud, Laboratoire de Math\'ematiques, b\^atiment 425\\
91405 Orsay Cedex, France}\\
{\it E-mail:} rcerf@math.u-psud.fr\\
\vskip0.5cm and\\
\vskip0.5cm {\Large Marie Th\'eret}\\
{\it \'Ecole Normale Sup\'erieure, D\'epartement Math\'ematiques et
Applications, 45 rue d'Ulm\\ 75230 Paris Cedex 05, France}\\
{\it E-mail:} marie.theret@ens.fr
\end{center}

\noindent
{\bf Abstract:}
We consider the standard first passage percolation model in the rescaled graph
$\ZZ^d/n$ for $d\geq 2$, and a domain $\O$ of boundary $\G$ in
$\RR^d$. Let $\G^1$ and $\G^2$ be two disjoint open subsets of $\G$, representing
the parts of $\G$ through which some water can enter and escape
from $\O$. We investigate the asymptotic behaviour of the flow $\phi_n$
through a discrete version $\O_n$ of $\O$ between the corresponding discrete
sets $\G^1_n$ and $\G^2_n$. We prove that under some
conditions on the regularity of the domain and on the law of the capacity of
the edges, the lower large deviations of $\phi_n/ n^{d-1}$ below a
certain constant are of surface order.\\



%%%%%%%%%%%%%%%%%%%%%%%%%%%%%%%%%%%%%%%%%%%%%%%%%%%%%%%%%%%%%%%%%%%%%%%%%%

\section{First definitions and main result}

We use many notations introduced in \cite{Kesten:StFlour} and
\cite{Kesten:flows}. Let $d\geq2$. We consider the graph $(\mathbb{Z}^{d}_n,
\mathbb E ^{d}_n)$ having for vertices $\mathbb Z ^{d}_n = \ZZ^d/n$ and for edges
$\mathbb E ^{d}_n$, the set of pairs of nearest neighbours for the standard
$L^{1}$ norm. With each edge $e$ in $\mathbb{E}^{d}_n$ we associate a random
variable $t(e)$ with values in $\mathbb{R}^{+}$. We suppose that the family
$(t(e), e \in \mathbb{E}^{d}_n)$ is independent and identically distributed,
with a common law $\Lambda$: this is the standard model of
first passage percolation on the graph $(\mathbb{Z}^d_n,
\mathbb{E}^d_n)$. We interpret $t(e)$ as the capacity of the edge $e$; it
means that $t(e)$ is the maximal amount of fluid that can go through the
edge $e$ per unit of time.

\begin{figure}[!ht]
\centering
\begin{picture}(0,0)%
\includegraphics{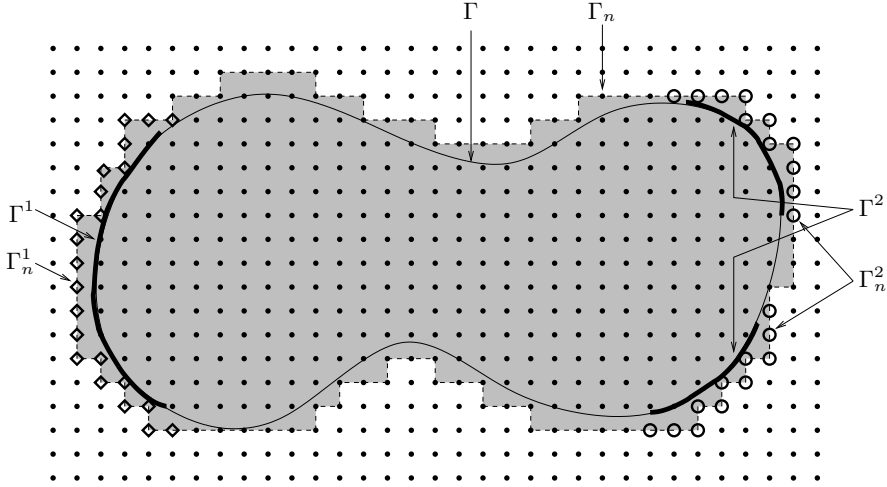}%
\end{picture}%
\setlength{\unitlength}{1973sp}%
\begingroup\makeatletter\ifx\SetFigFont\undefined%
\gdef\SetFigFont#1#2#3#4#5{%
  \reset@font\fontsize{#1}{#2pt}%
  \fontfamily{#3}\fontseries{#4}\fontshape{#5}%
  \selectfont}%
\fi\endgroup%
\begin{picture}(10380,6039)(1261,-6694)
\put(11626,-3361){\makebox(0,0)[lb]{\smash{{\SetFigFont{9}{10.8}{\rmdefault}{\mddefault}{\updefault}{\color[rgb]{0,0,0}$\Gamma^2$}%
}}}}
\put(1276,-3436){\makebox(0,0)[rb]{\smash{{\SetFigFont{9}{10.8}{\rmdefault}{\mddefault}{\updefault}{\color[rgb]{0,0,0}$\Gamma^1$}%
}}}}
\put(1276,-4036){\makebox(0,0)[rb]{\smash{{\SetFigFont{9}{10.8}{\rmdefault}{\mddefault}{\updefault}{\color[rgb]{0,0,0}$\Gamma^1_n$}%
}}}}
\put(11626,-4261){\makebox(0,0)[lb]{\smash{{\SetFigFont{9}{10.8}{\rmdefault}{\mddefault}{\updefault}{\color[rgb]{0,0,0}$\Gamma^2_n$}%
}}}}
\put(6751,-886){\makebox(0,0)[b]{\smash{{\SetFigFont{9}{10.8}{\rmdefault}{\mddefault}{\updefault}{\color[rgb]{0,0,0}$\Gamma$}%
}}}}
\put(8401,-886){\makebox(0,0)[b]{\smash{{\SetFigFont{9}{10.8}{\rmdefault}{\mddefault}{\updefault}{\color[rgb]{0,0,0}$\Gamma_n$}%
}}}}
\end{picture}%
\caption{Domain $\Omega$.}
\label{chapitre7domaine2}
\end{figure}

We consider an open bounded connected subset $\O$ of $\RR^d$ such that
the boundary $\G = \p \O$ of $\O$ is piecewise
of class $\C^1$ (in particular $\G$ has finite area: $\H^{d-1}(\G)
<\infty$). It means that $\G$ is included in the union of a finite number of
hypersurfaces of class $\C^1$, i.e., in the union of a finite number of
$C^1$ submanifolds of~$\RR^d$ of codimension~$1$. Let $\G^1$, $\G^2$ be two disjoint subsets of $\G$ that are
open in $\G$.
We want to define the maximal flow from $\G^1$ to $\G^2$ through $\O$ for
the capacities $(t(e), e\in \EE^d_n)$. We consider a discrete version
$(\O_n, \G_n, \G^1_n, \G^2_n)$ of $(\O, \G, \G^1,\G^2)$ defined by:
$$ \left\{ \begin{array}{l} \O_n \,=\, \{ x\in \ZZ^d_n \,|\,
    d_{\infty}(x,\O) <1/n  \}\,,\\ \G_n \,=\, \{ x\in
    \O_n \,|\, \exists y \notin \O_n\,,\,\, \langle x,y \rangle \in \EE^d_n
    \}\,,\\ \G^i_n \,=\, \{ x\in \G_n \,|\, d_\infty (x, \G^i) <1/n \,,\,\,
    d_\infty (x, \G^{3-i}) \geq 1/n \} \textrm{ for }i=1,2 \,,  \end{array}  \right.  $$
where $d_{\infty}$ is the $L^{\infty}$-distance, the notation $\langle
x,y\rangle$ corresponds to the edge of endpoints $x$ and $y$ (see
figure \ref{chapitre7domaine2}).

We shall study the maximal flow from $\G^1_n$ to $\G^2_n$ in $\O_n$.
Let us define
properly  the maximal flow $\phi(F_1 \rightarrow F_2 \textrm{ in } C)$ from
$F_1$ to $F_2$ in $C$, for $C \subset \mathbb{R}^d$ (or by commodity the
corresponding graph $C\cap \mathbb{Z}^d/n$). We will say that an edge
$e=\langle x,y\rangle$ belongs to a subset $A$ of $\mathbb{R}^{d}$, which
we denote by $e\in A$, if the interior of the segment joining $x$ to $y$ is included in $A$. We define
$\widetilde{\mathbb{E}}_n^{d}$ as the set of all the oriented edges, i.e.,
an element $\widetilde{e}$ in $\widetilde{\mathbb{E}}_n^{d}$ is an ordered
pair of vertices which are nearest neighbours. We denote an element $\widetilde{e} \in \widetilde{\mathbb{E}}_n^{d}$ by $\langle \langle x,y \rangle \rangle$, where $x$, $y \in \mathbb{Z}_n^{d}$ are the endpoints of $\widetilde{e}$ and the edge is oriented from $x$ towards $y$. We consider the set $\mathcal{S}$ of all pairs of functions $(g,o)$, with $g:\mathbb{E}_n^{d} \rightarrow \mathbb{R}^{+}$ and $o:\mathbb{E}_n^{d} \rightarrow \widetilde{\mathbb{E}}_n^{d}$ such that $o(\langle x,y\rangle ) \in \{ \langle \langle x,y\rangle \rangle , \langle \langle y,x \rangle \rangle \}$, satisfying:
\begin{itemize}
\item for each edge $e$ in $C$ we have
$$0 \,\leq\, g(e) \,\leq\, t(e) \,,$$
\item for each vertex $v$ in $C \smallsetminus (F_1\cup F_2)$ we have
$$ \sum_{e\in C\,:\, o(e)=\langle\langle v,\cdot \rangle \rangle}
  g(e) \,=\, \sum_{e\in C\,:\, o(e)=\langle\langle \cdot ,v \rangle
    \rangle} g(e) \,, $$
\end{itemize}
where the notation $o(e) = \langle\langle v,. \rangle \rangle$ (respectively $o(e) = \langle\langle .,v \rangle \rangle$) means that there exists $y \in \mathbb{Z}_n^d$ such that $e = \langle v,y \rangle$ and $o(e) = \langle\langle v,y \rangle \rangle$ (respectively $o(e) = \langle\langle y,v \rangle \rangle$).
A couple $(g,o) \in \mathcal{S}$ is a possible stream in $C$ from
$F_1$ to $F_2$: $g(e)$ is the amount of fluid that goes through the edge $e$, and $o(e)$ gives the direction in which the fluid goes through $e$. The two conditions on $(g,o)$ express only the fact that the amount of fluid that can go through an edge is bounded by its capacity, and that there is no loss of fluid in the graph. With each possible stream we associate the corresponding flow
$$ \flow (g,o) \,=\, \sum_{ u \in F_2 \,,\,  v \notin C \,:\, \langle
  u,v\rangle \in \mathbb{E}_n^{d}} g(\langle u,v\rangle) \mathbbm{1}_{o(\langle u,v\rangle) = \langle\langle u,v \rangle\rangle} - g(\langle u,v\rangle) \mathbbm{1}_{o(\langle u,v\rangle) = \langle\langle v,u \rangle\rangle} \,. $$
This is the amount of fluid that crosses $C$ from $F_1$
  to $F_2$ if the fluid respects the stream $(g,o)$. The maximal flow through
  $C$ from $F_1$ to $F_2$ is the supremum of this quantity over all possible choices of streams
$$ \phi (F_1 \rightarrow F_2 \textrm{ in }C)  \,=\, \sup \{ \flow (g,o)\,|\,
  (g,o) \in \mathcal{S} \}  \,.$$

We recall that we consider an open bounded
connected subset $\O$ of $\RR^d$ whose boundary $\G$ is piecewise of
class $\C^1$, and two disjoint open subsets $\G_1$ and $\G^2$ of
$\G$. We denote by
$$ \phi_n \,=\, \phi (\G^1_n \rightarrow \G^2_n \textrm{ in } \O_n) $$
the maximal flow from $\G^1_n$ to $\G^2_n$ in $\O_n$. We will investigate
the asymptotic behaviour of $\phi_n/n^{d-1}$ when $n$ goes to infinity. More precisely, we will show that the lower large deviations of
$\phi_n/n^{d-1}$ below a constant $\phi_{\O}$ are of surface
order. The description of $\phi_{\O}$ will be given in section
\ref{chapitre7deflimite}, and $p_c(d)$ is the critical parameter for
the bond percolation on $\ZZ^d$. Here we state the precise theorem:
\begin{thm}
\label{chapitre7devinf}
If the law $\Lambda$ of the capacity of an edge admits an exponential moment:
$$\exists \theta >0 \qquad \int_{\RR^+} e^{\theta x} d\Lambda (x) \,<\, +\infty \,, $$
and if $\Lambda (0) < 1-p_c(d)$, then there exists a finite constant
$\phi_{\O}$ such that for all $\lambda < \phi_{\O}$,
$$ \limsup_{n\rightarrow \infty} \frac{1}{n^{d-1}}\log \PP [ \phi_n \leq
\lambda n^{d-1}  ] \,<\,0\,.$$
\end{thm}

\begin{rem}
The lower large deviations we obtain are of the relevant
order. Indeed, if all the edges in a flat layer that separates $\G^1_n$ from $\G^2_n$ in $\O_n$ have abnormally small capacity,
then $\phi_n$ will be abnormally small. Since the cardinality of such a set
of edges is $D' n^{d-1}$ for a constant $D'$, the probability of this event
is of order $\exp -D n^{d-1}$ for a constant $D$.
\end{rem}

\begin{rem}
The condition $\Lambda(0) <1-p_c(d)$ is optimal. Indeed, Zhang proved
in \cite{Zhang} that in the particular case where $d=3$ and $\O$ is a straight
cube of bottom $\G^1$ and top $\G^2$, if $\Lambda$ admits an exponential moment
and $\Lambda(0)=1-p_c(d)$, then $\lim_{n\rightarrow \infty} \phi_n/n^{d-1} =
0$ a.s. The heuristic is the following: if $\Lambda(0) \geq 1-p_c(d)$, then
the edges of capacity strictly positive do not percolate, and therefore they cannot
convey a strictly positive amount of fluid through $\O$ when $n$ goes to
infinity. Kesten obtained the first results about maximal flows in
this model in \cite{Kesten:flows} under a stronger hypothesis on
$\Lambda(0)$. Zhang succeeded in
relaxing the constraint on $\Lambda$ in his remarkable article \cite{Zhang07}.
\end{rem}

\begin{rem}
In the two companion papers \cite{CerfTheret09geo} and \cite{CerfTheret09sup}, we prove in fact that
$\phi_{\O}$ is the almost sure limit of $\phi_n / n^{d-1}$ when
$n$ goes to infinity, and that the upper large deviations of $\phi_n /
n^{d-1}$ above $\phi_{\O}$ are of volume order.
\end{rem}

%%%%%%%%%%%%%%%%%%%%%%%%%%%%%%%%%%%%%%%%%%%%%%%%%%%%%%%%%%%%%%%%%%%%%%%%%%%%%%%%

\section{Computation of $\phi_{\O}$}
\label{chapitre7deflimite}

\subsection{Geometric notations}

\begin{figure}[!ht]
\centering
\begin{picture}(0,0)%
\includegraphics{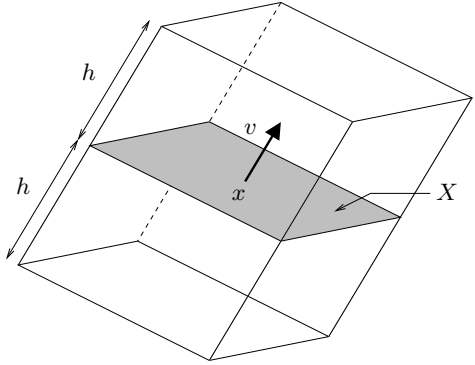}%
\end{picture}%
\setlength{\unitlength}{1973sp}%
\begingroup\makeatletter\ifx\SetFigFont\undefined%
\gdef\SetFigFont#1#2#3#4#5{%
  \reset@font\fontsize{#1}{#2pt}%
  \fontfamily{#3}\fontseries{#4}\fontshape{#5}%
  \selectfont}%
\fi\endgroup%
\begin{picture}(5874,4524)(2539,-6973)
\put(3676,-3436){\makebox(0,0)[rb]{\smash{{\SetFigFont{9}{10.8}{\rmdefault}{\mddefault}{\updefault}{\color[rgb]{0,0,0}$h$}%
}}}}
\put(2851,-4861){\makebox(0,0)[rb]{\smash{{\SetFigFont{9}{10.8}{\rmdefault}{\mddefault}{\updefault}{\color[rgb]{0,0,0}$h$}%
}}}}
\put(5701,-4111){\makebox(0,0)[rb]{\smash{{\SetFigFont{9}{10.8}{\rmdefault}{\mddefault}{\updefault}{\color[rgb]{0,0,0}$v$}%
}}}}
\put(5476,-4936){\makebox(0,0)[b]{\smash{{\SetFigFont{9}{10.8}{\rmdefault}{\mddefault}{\updefault}{\color[rgb]{0,0,0}$x$}%
}}}}
\put(7951,-4936){\makebox(0,0)[lb]{\smash{{\SetFigFont{9}{10.8}{\rmdefault}{\mddefault}{\updefault}{\color[rgb]{0,0,0}$X$}%
}}}}
\end{picture}%
\caption{Cylinder $\cyl(X,h)$.}
\label{chapitre7cylindrerect}
\end{figure}

We start with some geometric definitions. For a subset $X$ of
$\mathbb{R}^d$, we denote by $\mathcal{H}^s (X)$ the $s$-dimensional
Hausdorff measure of $X$ (we will use $s=d-1$ and $s=d-2$). The
$r$-neighbourhood $\V_i(X,r)$ of $X$ for the distance $d_i$, that can be
the Euclidean distance if $i=2$ or the $L^\infty$-distance if $i=\infty$,
is defined by
$$ \V_i (X,r) \,=\, \{ y\in \RR^d\,|\, d_i(y,X)<r\}\,.  $$
If $X$ is a subset of $\RR^d$ included in an hyperplane of $\RR^d$ and of
codimension $1$ (for example a non degenerate hyperrectangle), we denote by
$\hyp(X)$ the hyperplane spanned by $X$, and we denote by $\cyl(X, h)$ the
cylinder of basis $X$ and of height $2h$ defined by
$$ \cyl (X,h) \,=\, \{x+t v \,|\, x\in X \,,\,  t\in
[-h,h]    \}\,,$$
where $v$ is one of the two unit vectors orthogonal to $\hyp(X)$ (see
figure \ref{chapitre7cylindrerect}).

\begin{figure}[!ht]
\centering
\begin{picture}(0,0)%
\includegraphics{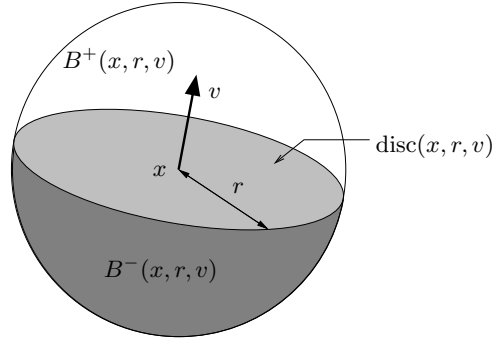}%
\end{picture}%
\setlength{\unitlength}{1973sp}%
\begingroup\makeatletter\ifx\SetFigFont\undefined%
\gdef\SetFigFont#1#2#3#4#5{%
  \reset@font\fontsize{#1}{#2pt}%
  \fontfamily{#3}\fontseries{#4}\fontshape{#5}%
  \selectfont}%
\fi\endgroup%
\begin{picture}(4598,4215)(3293,-6069)
\put(5776,-3061){\makebox(0,0)[lb]{\smash{{\SetFigFont{9}{10.8}{\rmdefault}{\mddefault}{\updefault}{\color[rgb]{0,0,0}$v$}%
}}}}
\put(4651,-2686){\makebox(0,0)[b]{\smash{{\SetFigFont{9}{10.8}{\rmdefault}{\mddefault}{\updefault}{\color[rgb]{0,0,0}$B^+(x,r,v)$}%
}}}}
\put(7876,-3736){\makebox(0,0)[lb]{\smash{{\SetFigFont{9}{10.8}{\rmdefault}{\mddefault}{\updefault}{\color[rgb]{0,0,0}$\disc(x,r,v)$}%
}}}}
\put(6076,-4261){\makebox(0,0)[lb]{\smash{{\SetFigFont{9}{10.8}{\rmdefault}{\mddefault}{\updefault}{\color[rgb]{0,0,0}$r$}%
}}}}
\put(5251,-4036){\makebox(0,0)[rb]{\smash{{\SetFigFont{9}{10.8}{\rmdefault}{\mddefault}{\updefault}{\color[rgb]{0,0,0}$x$}%
}}}}
\put(5176,-5311){\makebox(0,0)[b]{\smash{{\SetFigFont{9}{10.8}{\rmdefault}{\mddefault}{\updefault}{\color[rgb]{0,0,0}$B^-(x,r,v)$}%
}}}}
\end{picture}%
\caption{Ball $B(x,r)$.}
\label{chapitre7boule}
\end{figure}

For $x\in
\RR^d$, $r\geq 0$ and a unit vector $v$, we denote by $B(x,r)$ the closed
ball centered at $x$ of radius $r$, by $\disc (x,r,v)$ the
closed disc centered at $x$ of radius $r$ and normal vector $v$, and by
$B^+ (x,r,v)$ (respectively $B^- (x,r,v)$) the upper (respectively lower)
half part of $B(x,r)$ where the direction is determined by $v$ (see figure
\ref{chapitre7boule}), i.e.,
$$ B^+ (x,r,v) \,=\, \{y\in B(x,r) \,|\, (y-x)\cdot v \geq 0  \} \,,$$
$$ B^- (x,r,v) \,=\, \{y\in B(x,r) \,|\, (y-x)\cdot v \leq 0  \} \,.$$
We denote by $\alpha_d$ the volume of a unit ball in $\RR^d$, and
$\alpha_{d-1}$ the $\H^{d-1}$ measure of a unit disc.

%%%%%%%%%%%%%%%%%%%%%%%%%%%%%

\subsection{Flow in a cylinder}

Here are some particular definitions of flows through a box. It is important
to know them, because all our work consists in comparing the maximal flow
$\phi_n$ in $\O_n$ with the maximal flows in small cylinders.
Let $A$ be a non degenerate hyperrectangle,
i.e., a box of dimension $d-1$ in $\mathbb{R}^d$. All hyperrectangles will be
supposed to be closed in $\mathbb{R}^d$. We denote by
$v$ one of 
the two unit vectors orthogonal to $\hyp (A)$. For $h$ a
positive real number, we consider the cylinder $\cyl(A,h)$.
The set $\cyl(A,h) \smallsetminus \hyp (A)$ has two connected
components, which we denote by $\mathcal{C}_1(A,h)$ and
$\mathcal{C}_2(A,h)$. For $i=1,2$, let $A_i^h$ be
the set of the points in $\mathcal{C}_i(A,h) \cap \mathbb{Z}_n^d$ which have
a nearest neighbour in $\mathbb{Z}_n^d \smallsetminus \cyl(A,h)$:
$$ A_i^h\,=\,\{x\in \mathcal{C}_i(A,h) \cap
\mathbb{Z}_n^d \,|\, \exists y \in \mathbb{Z}_n^d \smallsetminus \cyl(A,h)
\,,\, \langle x,y \rangle \in \EE^d_n\}\,.$$
Let $T(A,h)$ (respectively $B(A,h)$) be the top
(respectively the bottom) of $\cyl(A,h)$, i.e.,
$$ T(A,h) \,=\, \{ x\in \cyl(A,h) \,|\, \exists y\notin \cyl(A,h)\,,\,\,
\langle x,y\rangle \in \mathbb{E}_n^d \textrm{ and }\langle x,y\rangle
\textrm{ intersects } A+hv  \}  $$
and
$$  B(A,h) \,=\, \{ x\in \cyl(A,h) \,|\, \exists y\notin \cyl(A,h)\,,\,\,
\langle x,y\rangle \in \mathbb{E}_n^d \textrm{ and } \langle x,y\rangle
\textrm{ intersects } A-hv  \} \,.$$
For a given realisation $(t(e),e\in \mathbb{E}_n^{d})$ we define the variable
$\tau (A,h) = \tau(\cyl(A,h), v)$ by
$$ \tau (A,h) \,=\,  \tau(\cyl(A,h), v)\,=\, \phi (A_1^h \rightarrow A_2^h
\textrm{ in } \cyl(A,h)) \,,$$
and the variable $\phi(A,h)= \phi(\cyl(A,h), v)$ by
$$ \phi(A,h) \,=\,\phi(\cyl(A,h), v) \,=\, \phi (B(A,h) \rightarrow T(A,h)
\textrm{ in }   \cyl(A,h))\,, $$ 
where $\phi(F_1 \rightarrow F_2 \textrm{ in } C)$ is the maximal
flow from $F_1$ to $F_2$ in $C$, for $C \subset \mathbb{R}^d$ (or by
commodity the corresponding graph $C\cap \mathbb{Z}^d/n$) defined
previously. The dependence in $n$ is implicit here, in fact we can
also write $\tau_n (A,h)$ and $\phi_n(A,h)$ if we want to emphasize
this dependence on the mesh of the graph.

%%%%%%%%%%%%%%%%%%%%%%%%%%%%%%%%%%

\subsection{Max-flow min-cut theorem}

The maximal flow
$\phi (F_1\rightarrow F_2 \textrm{ in } C)$ can be expressed differently
thanks to the max-flow min-cut theorem (see \cite{Bollobas}). We need some
definitions to state this result.
A path on the graph $\mathbb{Z}_n^{d}$ from $v_{0}$ to $v_{m}$ is a sequence $(v_{0}, e_{1}, v_{1},..., e_{m}, v_{m})$ of vertices $v_{0},..., v_{m}$ alternating with edges $e_{1},..., e_{m}$ such that $v_{i-1}$ and $v_{i}$ are neighbours in the graph, joined by the edge $e_{i}$, for $i$ in $\{1,..., m\}$.
A set $E$ of edges in $C$ is said to cut $F_1$ from $F_2$ in
$C$ if there is no path from $F_1$ to $F_2$ in $C \smallsetminus
E$. We call $E$ an $(F_1,F_2)$-cut if $E$ cuts $F_1$ from $F_2$ in $C$
and if no proper subset of $E$ does. With each set $E$ of edges we
associate its capacity which is the variable
$$ V(E)\, = \, \sum_{e\in E} t(e) \, .$$
The max-flow min-cut theorem states that
$$ \phi(F_1\rightarrow F_2 \textrm{ in } C) \, = \, \min \{ \, V(E) \, | \, E
\textrm{ is a } (F_1,F_2)\textrm{-cut} \, \} \, .$$

%%%%%%%%%%%%%%%%%%%%%%%%%%%%%%%%%%%%

\subsection{Definition of $\nu$}

The asymptotic behaviour of the rescaled expectation of $\tau_n (A,h)$
for large $n$ is well known, thanks to the almost subadditivity of
this variable. We recall the following result:
\begin{thm}
We suppose that
$$\int_{[0,+\infty[} x \, d\Lambda (x) \,<\,\infty \,.$$
Then for each unit vector $v$ there exists a constant $\nu (d,
\Lambda, v) = \nu(v)$ (the dependence on $d$ and $\Lambda$ is
implicit) such that for every non degenerate hyperrectangle $A$
orthogonal to $v$ and for every strictly positive
constant $h$, we have 
$$ \lim_{n\rightarrow \infty} \frac{\EE [\tau_n(A,h)]}{n^{d-1}
  \H^{d-1}(A)} \,=\, \nu(v) \,.$$
\end{thm}
For a proof of this proposition, see \cite{RossignolTheret08b}. We
emphasize the fact that the limit depends on the direction of $v$, but
not on $h$ nor on the hyperrectangle $A$ itself.

In fact, Rossignol and Th\'eret proved in \cite{RossignolTheret08b}
that under some moment conditions and/or some condition on $A$,
$\nu(v)$ is the limit of the rescaled
variable $\tau_n(A,h)/(n^{d-1} \H^{d-1}(A))$ almost surely and in $L^1$. We also know, thanks to
the works of Kesten \cite{Kesten:flows}, Zhang \cite{Zhang07} and
Rossignol and Th\'eret \cite{RossignolTheret08b} that the variable
$\phi_n(A,h)/ (n^{d-1}\H^{d-1}(A))$ satisfies the same law of large
numbers in the particular case where $A$ is a straight hyperrectangle,
i.e., a hyperrectangle of the form $\prod_{i=1}^{d-1}[0,k_i]\times
\{0\}$ for some $k_i>0$. In his article \cite{Zhang07}, Zhang obtains
a control on the number of edges in a minimal cutset. We will present
and use this result in section \ref{chapitre7seczhang}.

We recall some geometric properties of the map $\nu: v\in S^{d-1}
\mapsto \nu(v)$, under the only condition on $\Lambda$ that $\EE(t(e))<\infty$. They
have been stated in section 4.4 of \cite{RossignolTheret08b}. There exists a unit vector $v_0$ such that
$\nu(v_0)=0$ if and only if for all
 unit vector $v$, $\nu(v)=0$, and it happens if and only if $\Lambda(0) \geq
 1-p_c(d)$. This property has been proved by Zhang in \cite{Zhang}. Moreover, $\nu$ satisfies the weak
 triangle inequality, i.e., if $(ABC)$ is a non degenerate triangle in
 $\mathbb{R}^d$ and $v_A$, $v_B$ and $v_C$ are the
 exterior normal unit vectors to the sides $[BC]$, $[AC]$, $[AB]$ in the
 plane spanned by $A$, $B$, $C$, then
$$ \mathcal{H}^1 ([AB]) \nu(v_C) \,\leq\, \mathcal{H}^1 ([AC])
\nu(v_B) + \mathcal{H}^1 ([BC]) \nu(v_A) \,. $$
This implies that the homogeneous extension $\nu_0$ of $\nu$ to $\RR^d$,
defined by $ \nu_0(0) =0 $ and for all $ w$ in $\RR^d$, 
$$ \nu_0(w)\, =\, |w|_2 \nu (w/|w|_2)\,, $$
is a convex function; in particular, since $\nu_0$ is finite, it is
continuous on $\RR^d$. We denote by $\nu_{\min}$ (respectively
$\nu_{\max}$) the infimum (respectively supremum) of $\nu$ on
$S^{d-1}$.

The last result we recall is Theorem 3.9 in \cite{RossignolTheret08b}
concerning the lower large deviations of the variable $\tau_n(A,h)$
below $\nu(v)$:
\begin{thm}[Rossignol and Th\'eret]
\label{thmdevinftau}
We suppose that $\int_{[0,+\infty[} x \, d\Lambda (x) <\infty$ and
    that $\Lambda(0) <1-p_c(d)$. Then for every $\eps$ there exists a
    positive constant $K(d,\Lambda, \eps)$ such that for every unit
    vector $v$ and every non degenerate hyperrectangle $A$ orthogonal
    to $v$, there exists a constant $K'(d, \Lambda,
    A,\eps)$ such that for every strictly positive constant $h$ we have
$$ \PP \left[\frac{\tau_n (A,h)}{n^{d-1}\H^{d-1}(A)} \leq \nu(v) - \eps
      \right] \,\leq\,K'(d, \Lambda,
    A,\eps) \exp \left(- K(d,\Lambda, \eps) n^{d-1}\H^{d-1}(A)  \right) \,.$$
\end{thm}
We shall rely on this result for proving Theorem
\ref{chapitre7devinf}. Moreover, Theorem \ref{chapitre7devinf} is a
generalisation of Theorem \ref{thmdevinftau}, where we work in the domain
$\O$ instead of a parallelepiped.

%We consider the rescaled graph $(\ZZ^d_n, \EE^d_n)$ again, and for the rest
%of the article.

%%%%%%%%%%%%%%%%%%%%%%%%%%%%%%%%%%%%

\subsection{Definition of $\phi_{\O}$}

We give here a definition of $\phi_{\Omega}$ in
terms of the map $\nu$.
For a subset $F$ of $\RR^d$, we define the perimeter of $F$ in
$\O$ by
$$  \P(F, \O)  \,=\, \sup \left\{\int_F \di f(x) d \L^d(x), \, f\in \C_c^\infty
(\O, B(0,1))  \right\} \,, $$
where $\C_c^\infty (\O, B(0,1))$ is the set of the functions of class
$\C^{\infty}$ from $\RR^d$ to $B(0,1)$, the ball centered at $0$ and of
radius $1$ in $\RR^d$, having a compact support included in $\O$, and $\di$
is the usual divergence operator. The perimeter $\P(F)$ of $F$ is defined as $\P(F, \RR^d)$.
We denote by $\p F$ the boundary of $F$, and by
$\p^* F$ the reduced boundary of $F$. At any point $x$ of $\p^* F$, the set
$F$ admits a unit exterior normal vector $v_F(x)$ at $x$ in a measure
theoretic sense (for
definitions see for
example \cite{Cerf:StFlour}, section 13). For all $F \subset \RR^d$ of
finite perimeter in $\O$, we define
\begin{align*}
\I_{\O}(F) &\,=\, \int_{\p^* F \cap \O} \nu(v_F(x)) d\H^{d-1}(x) +
\int_{\G^2 \cap \p^* (F\cap \O)} \nu(v_{(F\cap \O)}(x)) d\H^{d-1}(x)\\ & \qquad \qquad +
\int_{\G^1 \cap \p^* (\O\smallsetminus F)} \nu(v_{\O}(x)) d\H^{d-1}(x)\,.
\end{align*}
If $\P(F, \O) = +\infty$, we define $\I_{\O}(F) =+\infty$. Finally, we define
$$ \phi_{\O} \,=\, \inf \{ \I_{\O}(F)\,|\, F\subset \RR^d  \} \,=\, \inf \{
\I_{\O}(F)\,|\, F\subset \O  \}\,. $$

\begin{figure}[!ht]
\centering
\begin{picture}(0,0)%
\includegraphics{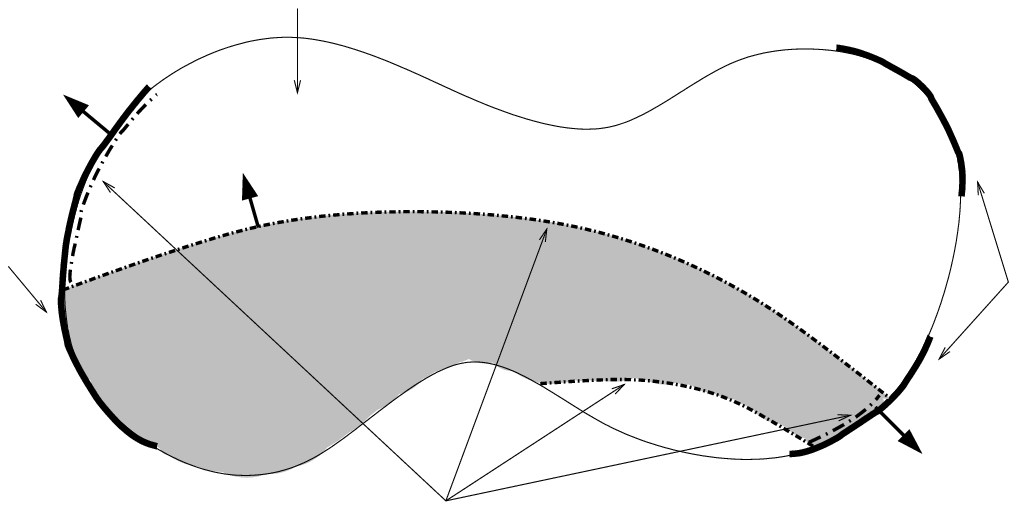}%
\end{picture}%
\setlength{\unitlength}{1973sp}%
\begingroup\makeatletter\ifx\SetFigFont\undefined%
\gdef\SetFigFont#1#2#3#4#5{%
  \reset@font\fontsize{#1}{#2pt}%
  \fontfamily{#3}\fontseries{#4}\fontshape{#5}%
  \selectfont}%
\fi\endgroup%
\begin{picture}(9780,5502)(1411,-6682)
\put(11176,-4261){\makebox(0,0)[lb]{\smash{{\SetFigFont{9}{10.8}{\rmdefault}{\mddefault}{\updefault}{\color[rgb]{0,0,0}$\Gamma^2$}%
}}}}
\put(1426,-4111){\makebox(0,0)[rb]{\smash{{\SetFigFont{9}{10.8}{\rmdefault}{\mddefault}{\updefault}{\color[rgb]{0,0,0}$\Gamma^1$}%
}}}}
\put(4276,-1411){\makebox(0,0)[b]{\smash{{\SetFigFont{9}{10.8}{\rmdefault}{\mddefault}{\updefault}{\color[rgb]{0,0,0}$\Omega$}%
}}}}
\put(3901,-3136){\makebox(0,0)[lb]{\smash{{\SetFigFont{9}{10.8}{\rmdefault}{\mddefault}{\updefault}{\color[rgb]{0,0,0}$v_F(x)$}%
}}}}
\put(3751,-3961){\makebox(0,0)[lb]{\smash{{\SetFigFont{9}{10.8}{\rmdefault}{\mddefault}{\updefault}{\color[rgb]{0,0,0}$x$}%
}}}}
\put(5176,-4411){\makebox(0,0)[lb]{\smash{{\SetFigFont{9}{10.8}{\rmdefault}{\mddefault}{\updefault}{\color[rgb]{0,0,0}$F$}%
}}}}
\put(2026,-2236){\makebox(0,0)[rb]{\smash{{\SetFigFont{9}{10.8}{\rmdefault}{\mddefault}{\updefault}{\color[rgb]{0,0,0}$v_{\Omega}(z)$}%
}}}}
\put(2626,-2986){\makebox(0,0)[lb]{\smash{{\SetFigFont{9}{10.8}{\rmdefault}{\mddefault}{\updefault}{\color[rgb]{0,0,0}$z$}%
}}}}
\put(10351,-5911){\makebox(0,0)[lb]{\smash{{\SetFigFont{9}{10.8}{\rmdefault}{\mddefault}{\updefault}{\color[rgb]{0,0,0}$v_{(F\cap \O)}(y)$}%
}}}}
\put(10201,-5461){\makebox(0,0)[b]{\smash{{\SetFigFont{9}{10.8}{\rmdefault}{\mddefault}{\updefault}{\color[rgb]{0,0,0}$y$}%
}}}}
\put(5776,-6586){\makebox(0,0)[b]{\smash{{\SetFigFont{9}{10.8}{\rmdefault}{\mddefault}{\updefault}{\color[rgb]{0,0,0}$(\partial F \cap \Omega) \cup (\Gamma^2 \cap \partial(F \cap \Omega)) \cup (\Gamma^1 \cap \partial (\Omega \smallsetminus F))$}%
}}}}
\end{picture}%
\caption{The set $(\p F \cap \O) \cup (\G^2 \cap \p (F\cap \O)) \cup (\G^1
  \cap \p(\O\smallsetminus F)) $.}
\label{chapitre7covering2}
\end{figure}

In the case where $\partial F$ is $\C^1$, $\I_{\O}(F)$ has the simpler
following expression:
\begin{align*}
\I_{\O}(F) &\,=\, \int_{\p F \cap \O} \nu(v_F(x)) d\H^{d-1}(x) +
\int_{\G^2 \cap \p (F\cap \O)} \nu(v_{(F\cap \O)}(x)) d\H^{d-1}(x)\\ & \qquad \qquad +
\int_{\G^1 \cap \p(\O\smallsetminus F)} \nu(v_{\O}(x)) d\H^{d-1}(x)\,.
\end{align*}
The localization of the set along which the previous integrals are done is
illustrated in figure \ref{chapitre7covering2}.
Since $\nu(v)$ is the average amount of fluid that can cross a hypersurface
of area one in the direction $v$ per unit of time, it can be interpreted as
the capacity of a unitary hypersurface orthogonal to $v$. Thus $\I_{\O}(F)$
can be interpreted as the capacity of  $(\p F \cap \O) \cup (\G^2 \cap \p
(F\cap \O)) \cup (\G^1 \cap \p(\O\smallsetminus F)) $.

\section{Sketch of the proof}

We are studying the lower large deviations of $\phi_n/n^{d-1}$: they are
controlled by what
happens around a minimal cutset. First, we will use the estimate of the
number of edges in a minimal cutset made by Zhang in \cite{Zhang07} to
restrict the problem to cutsets having a number of edges at most $c n^{d-1}$
for a constant $c$; we can then conclude that the minimal cutset is "near"
the boundary of a subset $F$ of $\O$ belonging to a compact space. By
making an adequate covering of this space, we need only to deal with a
finite number of sets and their neighbourhoods. We will then cover the
boundary of such a set $F$ by balls of very
small radius, such that $\p F$ is "almost flat" in each ball; we will also
show that if $\phi_n$ is smaller than $\phi_{\O}(1-\eps) n^{d-1}$ for some
positive $\eps$, then some local event happens in each ball of the
covering of $\p F$ (this event will be denoted by $G(B,v_{F}(x))$ for the
ball $B$ centered at $x\in \p F$). After that, we will construct a link
between this local event in a ball and the fact that the maximal flow
through a cylinder (included in the ball) is abnormally small. The lower
large deviations for the maximal flow through a cylinder are already
known (see \cite{RossignolTheret08b}). Finally, we calibrate the constants
to get Theorem \ref{chapitre7devinf}.

This proof is largely inspired by the methods used to study the Wulff
crystal in Ising model in dimension $d\geq 3$ (see for example
\cite{Cerf:StFlour}).

%%%%%%%%%%%%%%%%%%%%%%%%%%%%%%%%%%%%%%%%%%%%%%%%%%%%%%%%%%%%%%%%%%%%%%%%%%%%%%%%

\section{Number of edges in a minimal cutset and compactness}
\label{chapitre7seczhang}

We consider a $(\G^1_n, \G^2_n)$-cut $\E_n$ in $\O_n$ of minimal capacity,
i.e., $\phi_n = V(\E_n)$, and of
minimal number of edges (if there are more than one such cutset, we select
one of them by a deterministic algorithm). According to Theorem $1$ in
\cite{Zhang07}, adapted to our case as said in Remark $2$ in
\cite{Zhang07}, we know that:
\begin{thm}[Zhang]
If the law of the capacity of the edges admits an exponential moment, and
if $\Lambda(0)<1-p_c(d)$, then there exist constants $\beta_0 = \beta_0
(\Lambda,d)$, $C_i = C_i(\Lambda,d)$ for $i=1,2$ and $N=N(\Lambda,
d,\O,\G, \G^1, \G^2 )$ such that for all $\beta \geq
\beta_0$, for all $n\geq N$, we have
$$ \PP [\card(\E_n) \geq \beta n^{d-1}] \,\leq\, C_1 \exp (-C_2 \beta n^{d-1}) \,. $$
\end{thm}

We will always consider such large $n\geq N$. Thus with high probability the $(\G^1_n,\G^2_n)$-cut $\E_n$ has not "too much"
edges. We want now to change a little bit our point of view in order to
work with a
subset of $\RR^d$ rather than the cutset $\E_n$. We define for each edge $e$
the variable $t'(e) = \ind{\{ e \notin \E_n\} }$, and the set $\widetilde{E}_n
\subset \ZZ^d_n$ by
$$ \widetilde{E}_n = \{x\in \O_n \,|\, x \textrm{ is in an open
  cluster connected to } \G_n^1 \textrm{ for the percolation process
}(t'(e))_{e\in \O_n}  \} . $$
Then the edge boundary $\p^e \widetilde{E}_n$ of $\widetilde{E}_n$, defined
by
$$ \p^e \widetilde{E}_n \,=\, \{e=\langle x,y \rangle \in \ZZ^d_n \cap \O_n
\,|\, x \in \widetilde{E}_n \textrm{ and } y\notin \widetilde{E}_n   \}
\,,$$
is exactly equal to $\E_n$. We consider now the "non discrete version"
$E_n$ of $\widetilde{E}_n$ defined by
$$ E_n \,=\,\{x\in \O \,|\, d_{\infty}(x,\widetilde{E}_n) \leq 1/(2n)  \} \,=\, \left( \widetilde{E}_n + [-1/(2n),1/(2n)]^d \right) \cap \O  \,.$$
For all $F\subset \RR^d$, we recall that the perimeter of $F$ in $\O$ is
defined by
$$ \P (F, \O) \,=\, \sup \left\{\int_F \di f(x) d \L^d(x), \, f\in \C_c^\infty
(\O, B(0,1))  \right\} \,.$$
We know that if $\card (\E_n) \leq \beta n^{d-1} $, then $\P (E_n, \O)\leq \beta$.

We define
$$ \C_\beta \,=\, \{ F\subset \O \,|\, \P(F, \O) \leq \beta  \} \,,$$
endowed with the topology $L^1$ associated to the distance $d(F, F') =
\L^d(F \triangle F')$, where $F\triangle F'$ is the symmetric difference
between these two sets. For this topology the set $\C_\beta$ is
compact. With every $F$ in $\C_\beta$ we associate a positive $\eps_F$,
that we will choose later. The collection of sets $\V(F, \eps_F), \, F\in
\C_{\beta}$, where
$\V(F, \eps_F)$ is the neighbourhood of $F$ of size $\eps_F$ for the
distance defined previously, covers $\C_\beta$ so we can extract a finite
covering: $\C_\beta \subset
\cup_{i=1...N} \V(F_i, \eps_{F_i})$. We then obtain that for a fixed
$\beta \geq \beta_0$, for all $\lambda$ we have
\begin{align*}
\PP [\phi_n \leq \lambda n^{d-1}] & \,\leq\, e^{-\beta n^{d-1}} + \PP
[V(\E_n) \leq \lambda n^{d-1} \textrm{ and } \P(E_n, \O) \leq \beta ]\\
& \,\leq\, e^{-\beta n^{d-1}} + \sum_{i=1}^N \PP [V(\E_n) \leq \lambda
n^{d-1} \textrm{ and } \L^d(E_n \triangle F_i) \leq \eps_i]\,.
\end{align*}
It remains to study
$$\PP[V(\E_n) \leq \lambda n^{d-1} \textrm{ and } \L^d(E_n\triangle F)\leq
\eps_F]$$
for a generic $F$ in $\C_\beta$ and the corresponding $\eps_F$.

%%%%%%%%%%%%%%%%%%%%%%%%%%%%%%%%%%%%%%%%%%%%%%%%%%%%%%%%%%%%%%%%%%%%%%%%%%%%%%%%

\section{Covering of  $\p F$ by balls}

\subsection{Geometric tools}

We recall an important result about the
Minkowski content of a subset of $\RR^d$ (see for example Appendix A in
\cite{Cerf-Pisztora}). Whenever
$E$ is a closed $(d-1)$-rectifiable subset of $\RR^d$ (i.e.,
there exists a Lipschitz function mapping some bounded subset of
$\RR^{d-1}$ onto $E$), the Minkowski content of $E$, defined by
$$ \lim_{r\rightarrow 0} \frac{1}{2r} \L^d(\V_2(E, r))\,,   $$
exists and is equal to $\H^{d-1}(E)$.

We will also use the Vitali covering theorem for $\H^{d-1}$. A collection of sets $\U$ is called a Vitali
class for a Borel set $E$ of $\RR^d$ if for each $x\in
E$ and $\delta >0$, there exists a set $U\in \U$ containing $x$ such that
$0<\diam U < \delta$, where $\diam U$ is the diameter of the set $U$. We
now recall the Vitali covering theorem for $\H^{d-1}$ (see for instance \cite{FAL}, Theorem 1.10):
\begin{thm}
\label{chapitre7vitali}
Let $E$ be a $\H^{d-1}$ measurable subset of $\RR^d$ and $\U$ be a Vitali
class of closed sets for $E$. Then we may select a (countable) disjoint
sequence $(U_i)_{i\in I}$ from $\U$ such that 
$$ \textrm{either } \sum_{i\in I} (\diam U_i)^{d-1} \,=\, +\infty \textrm{ or
} \H^{d-1} (E\smallsetminus \cup_{i\in I} U_i) \,=\, 0 \,.$$
If $\H^{d-1} (E) <\infty$, then given $\eps >0$, we may also require that
$$ \H^{d-1} (E) \,\leq\, \frac{\alpha_{d-1}}{2^{d-1}} \sum_{i\in I} (\diam
U_i)^{d-1} \,. $$
\end{thm}
We recall next the Besicovitch differentiation theorem in $\RR^d$ (see for
example \cite{ASQU}):
\begin{thm}
\label{chapitre7besicovitch}
Let $\mathfrak{M}$ be a finite positive Radon measure on $\RR^d$. For any
Borel function $f\in L^1(\mathfrak{M})$, the quotient
$$ \frac{1}{\mathfrak{M}(B(x,r))} \int_{B(x,r)} f(y) d\mathfrak{M}(y) $$
converges $\mathfrak{M}$-almost surely towards $f(x)$ as $r$ goes to $0$.
\end{thm}
We state a result of covering that we will use in our study of the
lower deviations of $\phi_n$:
\begin{lem}
\label{chapitre7covering}
Let $F $ be a subset of $\O$ of finite perimeter. For every positive
constants $\delta$ and $\eta$, there exists a finite family of closed
disjoint balls
$(B_i)_{i\in I\cup J\cup K} = (B(x_i, r_i), v_i)_{i \in I\cup J \cup K}$
such that (the vector $v_i$ defines $B_i^-$) 
$$ \begin{array}{l} \forall i \in I \,,\,\, x_i \in \p ^*F\cap \O
    \,,\,\, r_i \in ]0,1[ \,,\,\, B_i \subset \O \,,\,\, \L^d((F\cap
    B_i) \triangle B_i^-) \,\leq\, \delta \alpha_d r_i^d \,,\\
\forall i\in J \,,\,\, x_i \in \G^{1}\cap \p^*(\O \smallsetminus  F) \,,\,\, r_i
\in ]0,1[ \,,\,\, \p \O \cap B_i \subset
\G^1 \,,\,\, \L^d((B_i\cap \O) \triangle B_i^-) \,\leq \, \delta
\alpha_d r_i^d \,,\\
\forall i \in K \,,\,\, x_i \in \G^{2} \cap \p^* F\,,\,\, r_i
\in ]0,1[ \,, \,\,  \p \O \cap B_i \subset \G^2  \,,\,\,
\L^d((F\cap B_i) \triangle B_i^-) \,\leq \, \delta \alpha_d r_i^d \,,\\
\end{array}  $$
and finally
$$\left| \I_{\O}(F)  - \sum_{i \in I \cup K} \alpha_{d-1} r_i^{d-1}
  \nu(v_F(x_i)) - \sum_{i \in J} \alpha_{d-1} r_i^{d-1} \nu(v_{\O}(x_i))
\right| \,\leq\, \eta\,.$$
\end{lem}

We will prove Lemma \ref{chapitre7covering} with the help of Theorems \ref{chapitre7vitali} and
\ref{chapitre7besicovitch}, following the proof of Lemma 14.6 in
\cite{Cerf:StFlour}. First notice that for $F\subset \O$, we have
\begin{align*}
\I_{\O}(F) &\,=\, \int_{\p^* F \cap \O} \nu(v_F(x)) d\H^{d-1}(x) +
\int_{\G^2 \cap \p^* F} \nu(v_{F}(x)) d\H^{d-1}(x)\\ & \qquad \qquad +
\int_{\G^1 \cap \p^* (\O\smallsetminus F)} \nu(v_{\O}(x)) d\H^{d-1}(x)\,.
\end{align*}
For $E$ a set of finite perimeter, we denote by $||\nabla_{\chi_E}||$ the
measure defined by
$$ \forall A \textrm{ Borel set in }\RR^d \qquad ||\nabla_{\chi_E}|| (A)
\,=\, \H^{d-1}(A\cap \p^* E) \,. $$
We consider a subset $F$ of $\O$ of finite perimeter.
We recall that the function $\nu: S^{d-1}\rightarrow \RR^+$ is
continuous. The map $x \in \p^* F\cap \O \mapsto v_F(x)$ is $||\nabla_{\chi_F
}||$-measurable, so we can apply the Besicovitch differentiation theorem in
$\RR^d$ to the maps $x\in \p^* F \cap \O \mapsto \nu(v_F(x))$ and $x\in
\p^* F \cap \O \mapsto 1$ to obtain that for $\H^{d-1}$-almost all $x\in
\p^* F \cap \O $
$$ \lim_{r\rightarrow 0} \frac{1}{\alpha_{d-1} r^{d-1}} \H^{d-1}(B(x,r)
\cap \p^* F \cap \O) \,=\, 1\,,  $$
$$ \lim_{r\rightarrow 0} \frac{1}{\alpha_{d-1} r^{d-1}} \int_{B(x,r) \cap
  \p^* F \cap \O}  \nu(v_F(y)) d\H^{d-1}(y) \,=\, \nu(v_F(x)) \,.$$
We denote by $\R_1$ the set of the points of $\p^* F \cap \O$ where the two
preceding identities hold simultaneously, thus $\H^{d-1} ((\p^* F \cap \O)
\smallsetminus \R_1) =0$. Similarly, let $\R_2$ be the set
of the points $x$ belonging to $\G^2 \cap \p^* F$ such that
$$ \lim_{r\rightarrow 0} \frac{1}{\alpha_{d-1} r^{d-1}} \H^{d-1}(B(x,r)
\cap \G^2 \cap \p^* F ) \,=\, 1\,,  $$
$$ \lim_{r\rightarrow 0} \frac{1}{\alpha_{d-1} r^{d-1}} \int_{B(x,r) \cap
  \G^2 \cap \p^* F}  \nu(v_F(y)) d\H^{d-1}(y) \,=\, \nu(v_F(x)) \,.$$
We also know that $\H^{d-1} ((\G^2 \cap \p^* F) \smallsetminus \R_2)
=0$. Since the map $x\in \G^1 \cap \p^* (\O \smallsetminus F) \mapsto
v_{\O}(x)$ is $||\nabla_{\chi_{\O}}||$-measurable, the same arguments
imply that the set $\R_3$ of the points $x$ of $\G^1 \cap \p^* (\O
\smallsetminus F)$ such that
$$ \lim_{r\rightarrow 0} \frac{1}{\alpha_{d-1} r^{d-1}} \H^{d-1}(B(x,r)
\cap \G^1 \cap \p^*(\O \smallsetminus F) ) \,=\, 1\,,  $$
$$ \lim_{r\rightarrow 0} \frac{1}{\alpha_{d-1} r^{d-1}} \int_{B(x,r) \cap
  \G^1 \cap \p^* (\O \smallsetminus F)}  \nu(v_{\O}(y)) d\H^{d-1}(y) \,=\,
\nu(v_{\O}(x)) \,,$$
satisfies $\H^{d-1}(\G^1 \cap \p^* (\O \smallsetminus F) \smallsetminus
\R_3)=0$. Moreover, from the theory of sets of finite perimeter (see for
example section 13 in \cite{Cerf:StFlour}), we know that
$$ \left\{ \begin{array}{ll} \forall x \in
\p^* F\,,\, & \lim_{r\rightarrow 0} r^{-d} \L^d(F \triangle
    B^-(x,r,v_{F}(x)))\,=\,0 \,,\\  \forall x \in \p^*(\O \smallsetminus
    F)\,,\,&  \lim_{r\rightarrow 0} r^{-d} \L^d(\O
    \triangle B^-(x,r,v_{\O}(x)))\,=\,0\,.  \end{array} \right. $$
We fix two parameters $\eta >0$ and $\delta>0$. For all $x\in \R_1$, there
exists a positive $r(x, \eta, \delta)$ such that for all $r < r(x, \eta,
\delta)$ we have
$$ |\H^{d-1} (B(x,r) \cap \p^* F \cap \O) -
    \alpha_{d-1}r^{d-1} | \,\leq\, \eta \alpha_{d-1} r^{d-1}\,,  $$
$$ \left|
      \frac{1}{\alpha_{d-1}r^{d-1}} \int_{B(x,r)\cap \p^* F \cap \O}
      \nu(v_F(y)) d\H^{d-1} (y) - \nu(v_F(x)) \right|\,\leq\,
    \eta\,, $$
$$ \L^d((F\cap B(x,r) )\triangle B^-(x, r, v_F(x)))\,\leq\, \delta
\alpha_d r^d \quad \textrm{and} \quad B(x,r) \,\subset\, \O  \,.$$
For all $x$ in $\R_2$, there
exists a positive $r(x, \eta, \delta)$ such that for all $r < r(x, \eta,
\delta)$ we have
$$ |\H^{d-1} (B(x,r) \cap \G^2 \cap \p^* F) -
    \alpha_{d-1}r^{d-1} | \,\leq\, \eta \alpha_{d-1} r^{d-1}\,,  $$
$$ \left|
      \frac{1}{\alpha_{d-1}r^{d-1}} \int_{B(x,r)\cap \G^2 \cap \p^* F}
      \nu(v_F(y)) d\H^{d-1} (y) - \nu(v_F(x)) \right|\,\leq\,
    \eta\,, $$
$$ \L^d((F\cap B(x,r) )\triangle B^-(x, r, v_F(x)))\,\leq\, \delta
\alpha_d r^d \quad \textrm{and} \quad B(x,r) \cap \G \,\subset\, \G^2
\,.$$
For all $x$ in $\R_3$, there
exists a positive $r(x, \eta, \delta)$ such that for all $r < r(x, \eta,
\delta)$ we have
$$ |\H^{d-1} (B(x,r) \cap \G^1 \cap \p^*(\O\smallsetminus  F)) -
    \alpha_{d-1}r^{d-1} | \,\leq\, \eta \alpha_{d-1} r^{d-1}\,,  $$
$$ \left|
      \frac{1}{\alpha_{d-1}r^{d-1}} \int_{B(x,r)\cap \G^1 \cap \p^*(\O
        \smallsetminus F)}
      \nu(v_{\O}(y)) d\H^{d-1} (y) - \nu(v_{\O}(x)) \right|\,\leq\,
    \eta\,, $$
$$ \L^d((\O  \cap B(x,r) )\triangle B^-(x, r, v_F(x)))\,\leq\, \delta
\alpha_d r^d \quad \textrm{and} \quad B(x,r) \cap \G \,\subset\, \G^1  \,.$$
The family of balls 
$$(B(x,r), x\in \R_1 \cup \R_2 \cup \R_3, r<r(x, \eta, \delta))$$
is a Vitali relation for $\R_1 \cup \R_2 \cup \R_3$. By the Vitali
covering theorem for $\H^{d-1}$, we may select from this collection of
balls a finite or countable collection of disjoint balls $B(x_i, r_i), i\in
I_1$ such that either
$$ \H^{d-1} \left( (\R_1 \cup \R_2 \cup \R_3) \smallsetminus \bigcup_{i\in I_1}
B(x_i, r_i) \right) \,=\,0 $$
or
$$  \sum_{i\in I_1} r_i^{d-1} \,=\, \infty\,. $$
We know that $\O$ and $F$ have finite perimeter, and that 
$$ (\p^* F \cap \O) \cup (\G^2 \cap \p^* F) \cup (\G^1 \cap \p^*(\O
  \smallsetminus F))  \,\subset\, \G \cup \p^* F\,,$$
so
\begin{align*}
(1-\eta) \sum_{i\in I_1} \alpha_{d-1} r_i^{d-1}& \,\leq\, \H^{d-1}\left((\p^* F \cap \O) \cup (\G^2 \cap \p^* F) \cup (\G^1 \cap \p^*(\O
  \smallsetminus F)) \right) \\
& \,\leq\, \H^{d-1} (\G \cup \p^* F) \,<\,\infty \,, 
\end{align*}
thus the first case occurs in the Vitali covering theorem, so we may select
a finite subset $I_2$ of $I_1$ such that
$$ \H^{d-1} \left( (\R_1 \cup \R_2 \cup \R_3) \smallsetminus \bigcup_{i\in I_2}
B(x_i, r_i) \right) \,\leq\, \eta  \H^{d-1} (\R_1 \cup \R_2 \cup \R_3) \,. $$
We claim that the collection of balls $(B(x_i, r_i), i\in I_2)$ enjoys the
desired properties. We define the sets
$$I \,=\, \{ i\in I_2 \,|\, x_i \in \p^* F \cap \O \}\,,$$
$$J \,=\, \{ i\in I_2 \,|\, x_i \in \G^1 \cap \p^* (\O \smallsetminus F) \}\,,$$
$$K \,=\, \{ i\in I_2 \,|\, x_i \in \G^2 \cap \p^* F \}\,,$$
and $v_i = v_F(x_i)$ for $i\in I \cup K$ and $v_i = v_{\O}(x_i)$ for $i\in
J$. Finally, we only have to check that
$$ \left| \I_{\O}(F)  - \sum_{i \in I\cup K} \alpha_{d-1} r_i^{d-1}
  \nu(v_F(x_i))  - \sum_{i \in J} \alpha_{d-1} r_i^{d-1} \nu(v_{\O}(x_i))
  \right| \,\leq\, \eta\,. $$
We recall that $\nu_{\max}$ is the supremum of $\nu$ over $S^{d-1}$; we have
\begin{align*}
\Bigg| \I_{\O}(F)  - \sum_{i\in I  \cup K} & \alpha_{d-1} r_i^{d-1} \nu(v_F(x_i))
  - \sum_{i\in J} \alpha_{d-1} r_i^{d-1} \nu(v_{\O}(x_i)) \Bigg| \\ 
& \,\leq\, \left| \int_{\R_1} \nu(v_F(y)) d\H^{d-1}(y) - \sum_{i\in I}
  \alpha_{d-1} r_i^{d-1} \nu(v_F(x_i))   \right|\\
& \qquad  + \left| \int_{\R_2}
  \nu(v_F(y)) d\H^{d-1}(y) - \sum_{i\in K} \alpha_{d-1} r_i^{d-1}
  \nu(v_F(x_i))   \right|\\
& \qquad + \left| \int_{\R_3} \nu(v_{\O}(y)) d\H^{d-1}(y) - \sum_{i\in J}
  \alpha_{d-1} r_i^{d-1} \nu(v_{\O}(x_i))   \right|\\
&  \,\leq\, \int_{\R_1  \smallsetminus \cup_{i\in I} B(x_i, r_i) }
\nu(v_F(y)) d\H^{d-1}(y) \\
& \qquad + \sum_{i\in I} \left|\int_{\R_1 \cap B(x_i,r_i)}
  \nu (v_F(y)) d\H^{d-1}(y) - \alpha_{d-1} r_i^{d-1} \nu(v_F(x)) \right|\\
&  \qquad + \int_{\R_2  \smallsetminus \cup_{i\in K} B(x_i, r_i) }
\nu(v_F(y)) d\H^{d-1}(y) \\
& \qquad + \sum_{i\in K} \left|\int_{\R_2 \cap B(x_i,r_i)}
  \nu (v_F(y)) d\H^{d-1}(y) - \alpha_{d-1} r_i^{d-1} \nu(v_F(x)) \right|\\
&  \qquad + \int_{\R_3  \smallsetminus \cup_{i\in J} B(x_i, r_i) }
\nu(v_{\O}(y)) d\H^{d-1}(y) \\
& \qquad + \sum_{i\in J} \left|\int_{\R_3 \cap B(x_i,r_i)}
  \nu (v_{\O}(y)) d\H^{d-1}(y) - \alpha_{d-1} r_i^{d-1} \nu(v_{\O}(x)) \right|\\
& \,\leq\,\eta \H^{d-1} (\R_1 \cup \R_2 \cup \R_3) \nu_{\max} +  \eta \sum_{i\in
  I \cup J \cup K}\alpha_{d-1}r_i^{d-1}\\
& \,\leq\, \eta \H^{d-1} (\R_1 \cup \R_2 \cup \R_3) \nu_{\max} \,+\, 2\eta
\H^{d-1} (\R_1 \cup \R_2 \cup \R_3)\\
& \,\leq\, \eta (\nu_{\max} + 2) (\P(F, \O) + \P(\O))\,.
\end{align*}
Since $(\nu_{\max} + 2) (\P(F, \O) + \P(\O))$ does not depend on $\eta$, we
have the required estimate.

%%%%%%%%%%%%%%%%%%%%%%%%%%%%%%

\subsection{Definition of a local event}

We consider a set $F$ in $\C_\beta$, and a positive $\eps_F$ that we have
to choose adequately. Thanks to Lemma \ref{chapitre7covering}, we
know that for every positive fixed $\delta$ and $\eta$, there exists a
finite family of closed disjoint balls 
$(B_i)_{i\in I\cup J\cup K} = (B(x_i, r_i), v_i)_{i \in I\cup J \cup K}$
such that (the vector $v_i$ defines $B_i^-$) 
$$ \begin{array}{l} \forall i \in I \,,\,\, x_i \in \p ^*F\cap \O
    \,,\,\, r_i \in ]0,1[ \,,\,\, B_i \subset \O \,,\,\, \L^d((F\cap
    B_i) \triangle B_i^-) \,\leq\, \delta \alpha_d r_i^d \,,\\
\forall i\in J \,,\,\, x_i \in \G^{1}\cap \p^*(\O \smallsetminus  F) \,,\,\, r_i
\in ]0,1[ \,,\,\, \p \O \cap B_i \subset
\G^1 \,,\,\, \L^d((B_i\cap \O) \triangle B_i^-) \,\leq \, \delta
\alpha_d r_i^d \,,\\
\forall i \in K \,,\,\, x_i \in \G^{2} \cap \p^* F\,,\,\, r_i
\in ]0,1[ \,, \,\,  \p \O \cap B_i \subset \G^2  \,,\,\,
\L^d((F\cap B_i) \triangle B_i^-) \,\leq \, \delta \alpha_d r_i^d \,,\\
\end{array}  $$
and finally
$$\left| \I_{\O}(F)  - \sum_{i \in I \cup K} \alpha_{d-1} r_i^{d-1}
  \nu(v_F(x_i)) - \sum_{i \in J} \alpha_{d-1} r_i^{d-1} \nu(v_{\O}(x_i))
\right| \,\leq\, \eta\,.$$
It is obvious that $\phi_{\O}<\infty$ because 
$$\phi_{\O} \,\leq\, \I_{\O} (\O) \,=\, \int_{\G^2 \cap \p^* \O}
\nu(v_{\O}(x)) d\H^{d-1}(x) \,\leq\, \nu_{\max} \H^{d-1}(\G^2) \,<\,
\infty \,.$$ 
We suppose for the rest of the article that $\phi_{\O} >0$ otherwise we do not have to
study any lower large deviations. We consider $\lambda < \phi_{\O}$. There
exists a positive $s$ (we can choose it smaller than $1$) such that
$\lambda \leq \phi_{\O} (1-2s) \leq \I_{\O}(F) (1-2s)$. We choose
$$  \eta \,=\, \frac{s\I_{\O}(F)}{4} \,, $$
and then we obtain that
\begin{align*}
 \bigg| \I_{\O}(F) & - \sum_{i\in I \cup K} \alpha_{d-1} r_i^{d-1}
   \nu(v_F(x_i))  - \sum_{i\in J} \alpha_{d-1} r_i^{d-1} \nu(v_{\O}(x_i))
  \bigg| \\ & \,\leq\, \left( \sum_{i\in I \cup K} \alpha_{d-1} r_i^{d-1}
    \nu(v_F(x_i)) +  \sum_{i\in J} \alpha_{d-1} r_i^{d-1}
    \nu(v_{\O}(x_i)) \right) \frac{s}{2} \,,
\end{align*}
and that
$$ \lambda \,\leq\,  \left( \sum_{i\in I  \cup K} \alpha_{d-1} r_i^{d-1}
    \nu(v_F(x_i)) +  \sum_{i\in J} \alpha_{d-1} r_i^{d-1}
    \nu(v_{\O}(x_i)) \right) (1-s)  \,.$$
Since the $(B_i)_{i\in I \cup J \cup K}$ are disjoint, we also know that
$$ \phi_n \,\geq\, \sum_{i\in I \cup J \cup K} V(\E_n \cap B_i) \,. $$
Then
\begin{align*}
\PP [ V(\E_n)& \leq \lambda  n^{d-1} \textrm{ and } \L^d(E_n \triangle F)
\leq \eps_F]\\ & \,\leq\, \PP \left[ \begin{array}{rl} \sum_{i\in I \cup J \cup K} V(\E_n \cap  B_i) \leq (1-s)& n^{d-1} \big(\sum_{i\in I  \cup K} \alpha_{d-1}
  r_i^{d-1} \nu(v_F(x_i))\\ & + \sum_{i\in J} \alpha_{d-1}
  r_i^{d-1} \nu(v_{\O}(x_i)) \big) \\  \textrm{and } \L^d(E_n \triangle
  F) & \leq \,
  \varepsilon_F \end{array}\right]\,.  
\end{align*}
From now on we choose $\varepsilon_F$ to be  
$$\eps_F \,=\, \min_{i\in I\cup J\cup K} \alpha_{d} r_i^{d} \delta\,, $$
for a fixed $\delta$ that we will choose later. For all $i \in I$, we
then have
$$ \L^d((E_n \cap B_i)\triangle B_i^-) \,\leq\, \L^d((F\cap B_i)\triangle
B_i^-) + \L^d(E_n \triangle F) \,\leq\, 2 \delta \alpha_d r_i^d \,.$$
We want to evaluate $\card (((E_n\cap B_i) \triangle B_i^-) \cap
\ZZ^d_n)$. It is equivalent to evaluate
$$ n^d \L^d(((E_n\cap B_i) \triangle B_i^-)\cap \ZZ^d_n + [-1/2n,
1/2n]^{d})  \,.$$
By definition, for all $x \in E_n  \cap \ZZ^d_n = \widetilde{E}_n$, $x+
[-1/2n,1/2n]^d \subset E_n$, so
\begin{align*} 
((E_n\cap B_i)  & \triangle B_i^-)\cap \ZZ^d_n + [-1/2n, 1/2n]^{d}\\
& \,\subset\, ((E_n\cap B_i) \triangle B_i^-) \cup (\V_{\infty}(B_i,
1/n)\smallsetminus
B_i) \cup (\V_{\infty}(B_i^-, 1/n)\smallsetminus B_i^-)\\
& \,\subset \, ((E_n\cap B_i)\triangle B_i^-) \cup (\V_{2}(B_i,
2d/n)\smallsetminus B_i) \cup (\V_{2}(B_i^-, 2d/n)\smallsetminus B_i^-) \,.
\end{align*}
Since $\p B_i$
and $\p B_i^-$ are very regular, the result about the Minkowski content implies that
$$ \lim_{n\rightarrow \infty} \frac{n}{2d} \L^d(\V_{2}(B_i,
2d/n)\smallsetminus B_i) \,=\, \H^{d-1} (\p B_i)$$
and
$$ \lim_{n\rightarrow \infty} \frac{n}{2d} \L^d(\V_{2}(B_i^-, 2d/n)\smallsetminus
B_i^-) \,=\, \H^{d-1} (\p B_i^-) \,.$$
For $n$ large enough, we then obtain that
$$ \L^d(((E_n\cap B_i) \triangle B_i^-)\cap \ZZ^d_n  + [-1/2n, 1/2n]^{d})
\,\leq\, 2 \delta \alpha_d r_i^d + \frac{4d(\H^{d-1}(\p B_i) + \H^{d-1}(\p
  B_i^-))}{n} \,,$$
and then for all $n$ large enough
\begin{align*}
\card(((E_n\cap B_i) \triangle B_i^-)\cap \ZZ^d_n) & \,\leq\, 2 \delta
\alpha_d r_i^d n^d + 4d(\H^{d-1}(\p B_i) + \H^{d-1}(\p B_i^-)) n^{d-1}\\
& \,\leq\, 4 \delta \alpha_d r_i^d n^d \,.
\end{align*}
For $i\in K$, exactly the same arguments imply that
$$ \card(((E_n\cap B_i) \triangle B_i^-)\cap \ZZ^d_n) \,\leq\, 4 \delta
\alpha_d r_i^d n^d   $$ 
for $n$ large enough.

\begin{figure}[!ht]
\centering
\begin{picture}(0,0)%
\includegraphics{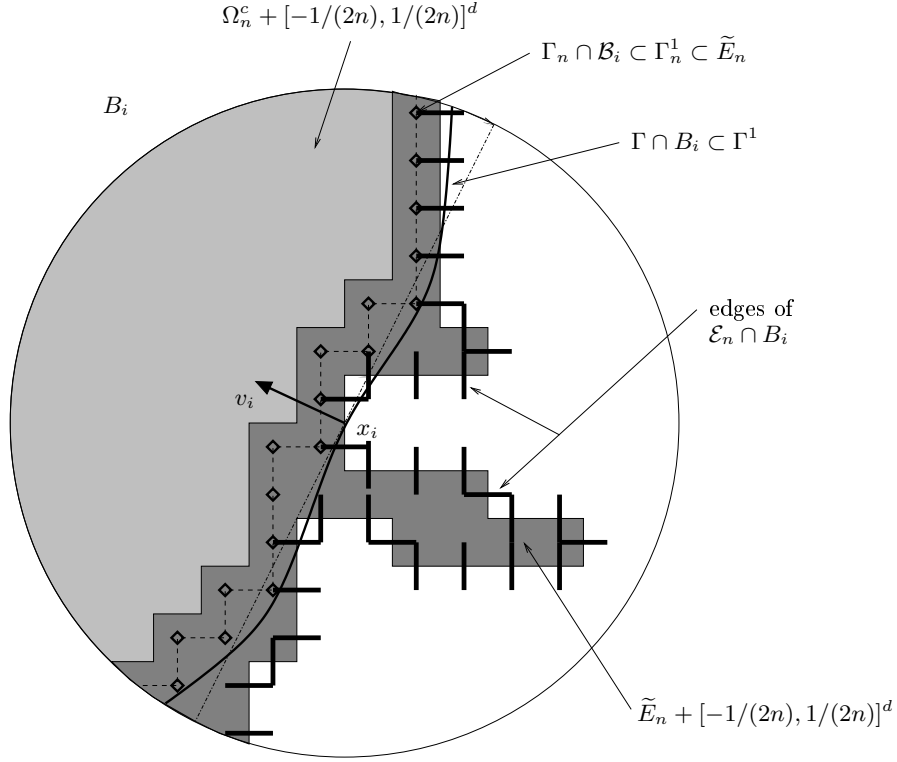}%
\end{picture}%
\setlength{\unitlength}{1973sp}%
\begingroup\makeatletter\ifx\SetFigFont\undefined%
\gdef\SetFigFont#1#2#3#4#5{%
  \reset@font\fontsize{#1}{#2pt}%
  \fontfamily{#3}\fontseries{#4}\fontshape{#5}%
  \selectfont}%
\fi\endgroup%
\begin{picture}(9739,9463)(1490,-9068)
\put(3001,-961){\makebox(0,0)[rb]{\smash{{\SetFigFont{9}{10.8}{\rmdefault}{\mddefault}{\updefault}{\color[rgb]{0,0,0}$B_i$}%
}}}}
\put(5851,-5011){\makebox(0,0)[lb]{\smash{{\SetFigFont{9}{10.8}{\rmdefault}{\mddefault}{\updefault}{\color[rgb]{0,0,0}$x_i$}%
}}}}
\put(10276,-3511){\makebox(0,0)[lb]{\smash{{\SetFigFont{9}{10.8}{\rmdefault}{\mddefault}{\updefault}{\color[rgb]{0,0,0}edges of}%
}}}}
\put(10276,-3841){\makebox(0,0)[lb]{\smash{{\SetFigFont{9}{10.8}{\rmdefault}{\mddefault}{\updefault}{\color[rgb]{0,0,0}$\E_n \cap B_i$}%
}}}}
\put(9301,-1411){\makebox(0,0)[lb]{\smash{{\SetFigFont{9}{10.8}{\rmdefault}{\mddefault}{\updefault}{\color[rgb]{0,0,0}$\G \cap B_i \subset \G^1$}%
}}}}
\put(8176,-286){\makebox(0,0)[lb]{\smash{{\SetFigFont{9}{10.8}{\rmdefault}{\mddefault}{\updefault}{\color[rgb]{0,0,0}$\G_n \cap \B_i \subset \G^1_n \subset \widetilde{E}_n$}%
}}}}
\put(9376,-8611){\makebox(0,0)[lb]{\smash{{\SetFigFont{9}{10.8}{\rmdefault}{\mddefault}{\updefault}{\color[rgb]{0,0,0}$\widetilde{E}_n + [-1/(2n), 1/(2n)]^d$}%
}}}}
\put(4576,-4636){\makebox(0,0)[rb]{\smash{{\SetFigFont{9}{10.8}{\rmdefault}{\mddefault}{\updefault}{\color[rgb]{0,0,0}$v_i$}%
}}}}
\put(5776,164){\makebox(0,0)[b]{\smash{{\SetFigFont{9}{10.8}{\rmdefault}{\mddefault}{\updefault}{\color[rgb]{0,0,0}$\O_n^c + [-1/(2n), 1/(2n)]^d $}%
}}}}
\end{picture}%
\caption{A ball $B_i$ for $i\in J$.}
\label{chapitre7Bi}
\end{figure}

We study now what happens in the balls $B_i$ for $i\in J$. We recall that
$\widetilde{E}_n = E_n \cap \ZZ^d_n$. We define
$\widetilde{E}_n' = \widetilde{E}_n \cup \O_n^c $ (where $\O_n^c = \ZZ^d_n
\smallsetminus \O_n$) and
$E_n'=\widetilde{E}_n' + [-1/(2n),1/(2n)]^{d-1}$. Then $E_n' \cap \O =
E_n$. In a ball $B_i$, we have $\p^e \widetilde{E}_n' \cap B_i =
\E_n \cap B_i$. Indeed, we know that  $\G \cap B_i \subset \G^1$. The sets
$\G^1 $ and $\G^2$ are open in $\G$ and disjoint, so $\G^1 \cap
\overline{\G^2} = \varnothing$, where $\overline{\G^2}$ is the adherence of
$\G^2$, and then $B_i \cap \overline{\G^2} =\varnothing$. Since $B_i$ is
closed, we obtain that $d(B_i, \overline{\G^2}) >0$, and thus for
$n$ large enough, $\G_n \cap B_i \subset \G^1_n$. Moreover, we know that
$\G^{1}_n \subset \widetilde{E}_n \subset \widetilde{E}_n'$. We obtain that
$\p^e \widetilde{E}_n' \cap \O_n^c \cap B_i = \varnothing$, i.e., all the
edges of $\p^e  \widetilde{E}_n'$ in $B_i$ have both endpoints in
$\O_n$ (see figure \ref{chapitre7Bi}).
Now we have
\begin{align*}
\L^d((E_n' \cap B_i) \triangle B_i^+) & \,\leq\, \L^d ( (E_n'\cap B_i)
\triangle (\O^c \cap B_i)  )+ \L^d((\O^c \cap B_i)\triangle B_i^+) \\
& \,\leq\, \L^d(E_n' \cap B_i \cap \O) + \L^d ((\O^c\smallsetminus
E_n')\cap B_i) +  \L^d((\O \cap B_i)\triangle B_i^-)\\
& \,\leq\,  \L^d(E_n\triangle F) + \L^d(\V_\infty (\G, 1/n)
\cap B_i) + \delta \alpha_{d} r_i^{d}\\
& \,\leq\, \eps_F + \L^d(\V_\infty (\G, 1/n)
\cap B_i) + \delta\alpha_{d} r_i^{d}\\
& \,\leq\, 3 \delta \alpha_{d} r_i^{d}\,,
\end{align*}
for $n$ large enough, where the last inequality is a consequence of the properties
of the Minkowski content. As previously, we
obtain that for $n$ large enough,
$$ \card(((E'_n\cap B_i) \triangle B_i^+)\cap \ZZ^d_n) \,\leq\, 4 \delta
\alpha_d r_i^d n^d  \,. $$
We conclude that for $n$ large enough,
\begin{align*}
& \PP [ V(\E_n) \leq  \lambda  n^{d-1}  \textrm{ and } \L^d(E_n \triangle F)
\leq \eps_F] \\
& \,\,\leq\, \sum_{i\in I} \PP \left[ \begin{array}{c} V(\p^e \widetilde{E}_n \cap B_i) \leq
(1-s)\alpha_{d-1} r_i^{d-1} \nu(v_F(x_i)) \textrm{ and}\\ \card(
(\widetilde{E}_n\cap B_i) \triangle (B_i^- \cap \ZZ^d_n)) \leq  4 \delta
\alpha_d r_i^d n^d \end{array} \right]\\ 
& \,\,\, + \sum_{i\in J} \PP \left[ \begin{array}{c} V(\p^e \widetilde{E}'_n \cap B_i) \leq
(1-s)\alpha_{d-1} r_i^{d-1} \nu(v_F(x_i)) \textrm{ and}\\\card( (
\widetilde{E}'_n\cap B_i) \triangle (B_i^+ \cap \ZZ^d_n)) \leq  4 \delta
\alpha_d r_i^d n^d \end{array} \right]\\
& \,\,\, + \sum_{i\in K}  \PP \left[ \begin{array}{c} V(\p^e \widetilde{E}_n \cap B_i) \leq
(1-s)\alpha_{d-1} r_i^{d-1} \nu(v_F(x_i)) \textrm{ and} \\ \card((\widetilde{E}_n\cap B_i)
\triangle (B_i^- \cap \ZZ^d_n)) \leq  4 \delta \alpha_d r_i^d n^d \end{array} \right]\\
&\leq\,  \sum_{i\in I\cup J \cup K} \PP [G(x_i,r_i,v_i )] \,,\\
\end{align*}
where $G(x,r, v)$ is the event that there exists a set $U \subset B\cap
\ZZ^d_n$ such that:
$$ \left\{ \begin{array}{l}
\card( U \triangle  B^- ) \,\leq\, 4\delta \alpha_d r^d  n^{d}\,,  \\
V(\p^e U \cap B) \,\leq\, (\alpha_{d-1} r^{d-1} \nu (v(x))) (1-s) n^{d-1}\,.
\end{array} \right.$$
Notice that this event depends only on the edges in $B = B(x,r)$. This event seems to be complicated, but indeed when $G(x,r,v)$ happens, it
means in a sense that the flow between the lower half part of $B(x,r)$ (for the
direction $v$) and the upper half part of $B$ is abnormally small. We
will examine the consequence of the event $G(x,r,v)$ over the maximal
flow in $B(x,r)$ in the next section.

%%%%%%%%%%%%%%%%%%%%%%%%%%%%%%%%%%%%%%%%%%%%%%%%%%%%%%%%%%%%%%%%%%%%%%%%%%%%%%%%

\section{Surgery in a ball to define an almost flat cutset}

\begin{figure}[!ht]
\centering
\begin{picture}(0,0)%
\includegraphics{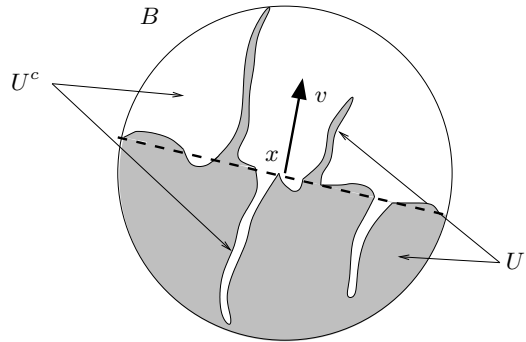}%
\end{picture}%
\setlength{\unitlength}{1973sp}%
\begingroup\makeatletter\ifx\SetFigFont\undefined%
\gdef\SetFigFont#1#2#3#4#5{%
  \reset@font\fontsize{#1}{#2pt}%
  \fontfamily{#3}\fontseries{#4}\fontshape{#5}%
  \selectfont}%
\fi\endgroup%
\begin{picture}(5880,4221)(2311,-6075)
\put(3826,-2086){\makebox(0,0)[rb]{\smash{{\SetFigFont{9}{10.8}{\rmdefault}{\mddefault}{\updefault}{\color[rgb]{0,0,0}$B$}%
}}}}
\put(8176,-5161){\makebox(0,0)[lb]{\smash{{\SetFigFont{9}{10.8}{\rmdefault}{\mddefault}{\updefault}{\color[rgb]{0,0,0}$U$}%
}}}}
\put(2326,-2911){\makebox(0,0)[rb]{\smash{{\SetFigFont{9}{10.8}{\rmdefault}{\mddefault}{\updefault}{\color[rgb]{0,0,0}$U^c$}%
}}}}
\put(5326,-3811){\makebox(0,0)[rb]{\smash{{\SetFigFont{9}{10.8}{\rmdefault}{\mddefault}{\updefault}{\color[rgb]{0,0,0}$x$}%
}}}}
\put(5776,-3061){\makebox(0,0)[lb]{\smash{{\SetFigFont{9}{10.8}{\rmdefault}{\mddefault}{\updefault}{\color[rgb]{0,0,0}$v$}%
}}}}
\end{picture}%
\caption{Event $G(x,r,v)$.}
\label{chapitre7filament2}
\end{figure}

We consider a fixed ball $B=B(x, r)$ and a unit vector $v$
(corresponding to one generic ball of the previous covering).
We want to interpret the event $G(x,r, v)$ in
term of the maximal flow through a cylinder whose basis is a disc, included in the
ball $B$, and oriented along the direction $v$. We define
$$ \gamma_{\max} \,=\, \rho r \,, $$
where $\rho$ is a constant depending on $\delta$ and $B$ which we can
imagine very small, it will be chosen later. The constant $\gamma_{\max}$
is in fact the height of
the cylinder we are constructing, namely
$$ \C \,=\, \cyl (\disc(x, r', v), \gamma_{\max}) \,.$$
We want $\C$ to be included in $B$, so we choose
$$ r' \,=\, r \cos (\arcsin \rho) \,. $$
We would like to analyse the implication of the event $G(x,r,v)$ on the
flow $\phi_{\C}$ between the top and the bottom of $\C$ for the
direction $v$ (we will define it properly soon). As we said previously, the event $G(x,r,v)$ means that the maximal flow between a set $U$
that "looks like" $B^-$ (for the direction given by $v$) and the set $U^c$
that "looks like" $B^+$ is a bit too small. Here "looks like" means that
$B^-$ and
$U$ are closed in volume, but the set $U$ might have some
thin strands (of small volume, but that can be long) that go deeply into
$B^+$ and symmetrically the set $U^c$ might have some thin strands that go
deeply into $B^-$ (see figure \ref{chapitre7filament2}).
What we have to do to control $\phi_{\C}$ is to cut these
strands: by adding edges to $\p^e U$ at a fixed height in $\C$ to close
the strands, we obtain a cutset in $\C$. The point is that we have to
control the capacity of these edges we have added to $\p^e U$. This is the
reason why we choose the
height at which we add edges to be sure we add not too many edges, and then
we control their capacity thanks to a property of independence.

\begin{figure}[!ht]
\centering
\begin{picture}(0,0)%
\includegraphics{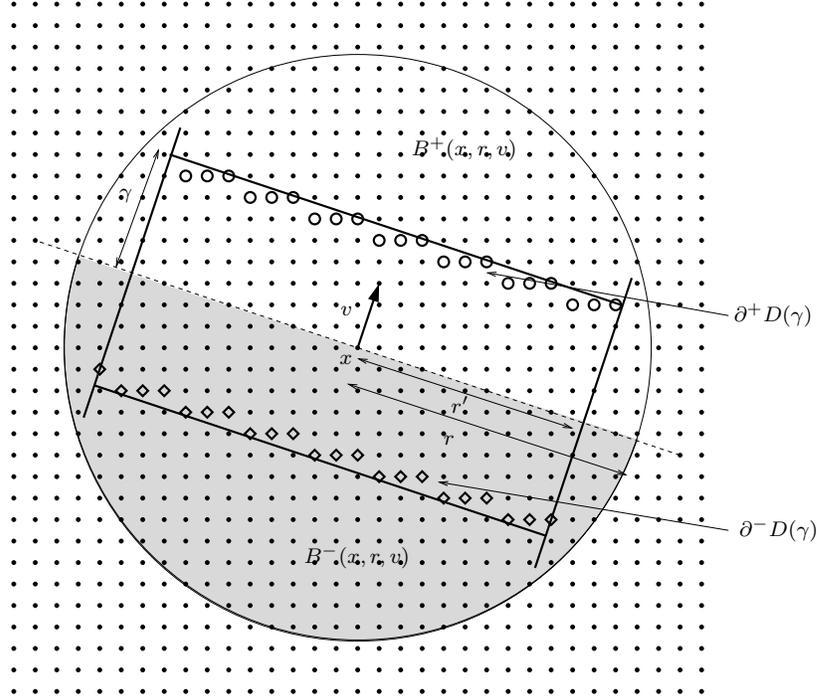}%
\end{picture}%
\setlength{\unitlength}{1776sp}%
\begingroup\makeatletter\ifx\SetFigFont\undefined%
\gdef\SetFigFont#1#2#3#4#5{%
  \reset@font\fontsize{#1}{#2pt}%
  \fontfamily{#3}\fontseries{#4}\fontshape{#5}%
  \selectfont}%
\fi\endgroup%
\begin{picture}(10174,9666)(1467,-11194)
\put(3151,-4261){\makebox(0,0)[rb]{\smash{{\SetFigFont{8}{9.6}{\rmdefault}{\mddefault}{\updefault}{\color[rgb]{0,0,0}$\gamma$}%
}}}}
\put(7051,-3661){\makebox(0,0)[lb]{\smash{{\SetFigFont{8}{9.6}{\rmdefault}{\mddefault}{\updefault}{\color[rgb]{0,0,0}$B^+(x,r,v)$}%
}}}}
\put(6226,-5911){\makebox(0,0)[rb]{\smash{{\SetFigFont{8}{9.6}{\rmdefault}{\mddefault}{\updefault}{\color[rgb]{0,0,0}$v$}%
}}}}
\put(6226,-6586){\makebox(0,0)[rb]{\smash{{\SetFigFont{8}{9.6}{\rmdefault}{\mddefault}{\updefault}{\color[rgb]{0,0,0}$x$}%
}}}}
\put(7576,-7711){\makebox(0,0)[b]{\smash{{\SetFigFont{8}{9.6}{\rmdefault}{\mddefault}{\updefault}{\color[rgb]{0,0,0}$r$}%
}}}}
\put(5551,-9361){\makebox(0,0)[lb]{\smash{{\SetFigFont{8}{9.6}{\rmdefault}{\mddefault}{\updefault}{\color[rgb]{0,0,0}$B^-(x,r,v)$}%
}}}}
\put(11551,-5986){\makebox(0,0)[lb]{\smash{{\SetFigFont{8}{9.6}{\rmdefault}{\mddefault}{\updefault}{\color[rgb]{0,0,0}$\partial^+ D(\gamma)$}%
}}}}
\put(11626,-8986){\makebox(0,0)[lb]{\smash{{\SetFigFont{8}{9.6}{\rmdefault}{\mddefault}{\updefault}{\color[rgb]{0,0,0}$\partial^- D(\gamma)$}%
}}}}
\put(7726,-7261){\makebox(0,0)[b]{\smash{{\SetFigFont{8}{9.6}{\rmdefault}{\mddefault}{\updefault}{\color[rgb]{0,0,0}$r'$}%
}}}}
\end{picture}%
\caption{Representation of $D(\gamma)$.}
\label{chapitre7boite}
\end{figure}

\begin{figure}[!ht]
\centering
\begin{picture}(0,0)%
\includegraphics{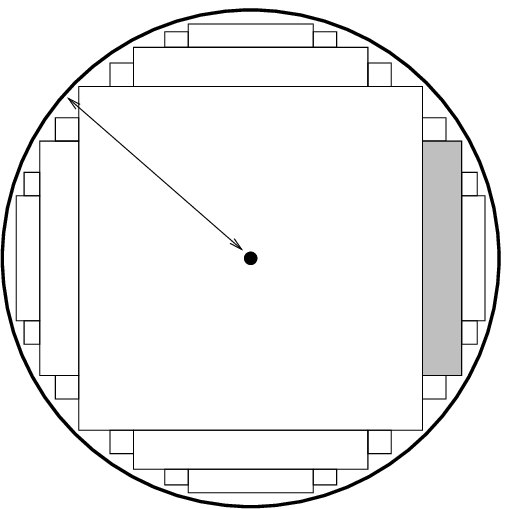}%
\end{picture}%
\setlength{\unitlength}{1973sp}%
\begingroup\makeatletter\ifx\SetFigFont\undefined%
\gdef\SetFigFont#1#2#3#4#5{%
  \reset@font\fontsize{#1}{#2pt}%
  \fontfamily{#3}\fontseries{#4}\fontshape{#5}%
  \selectfont}%
\fi\endgroup%
\begin{picture}(4814,4812)(3594,-6967)
\put(7801,-4336){\makebox(0,0)[b]{\smash{{\SetFigFont{9}{10.8}{\rmdefault}{\mddefault}{\updefault}{\color[rgb]{0,0,0}$A_i$}%
}}}}
\put(5926,-4861){\makebox(0,0)[rb]{\smash{{\SetFigFont{9}{10.8}{\rmdefault}{\mddefault}{\updefault}{\color[rgb]{0,0,0}$x$}%
}}}}
\put(5401,-3886){\makebox(0,0)[lb]{\smash{{\SetFigFont{9}{10.8}{\rmdefault}{\mddefault}{\updefault}{\color[rgb]{0,0,0}$r'$}%
}}}}
\put(3826,-2536){\makebox(0,0)[rb]{\smash{{\SetFigFont{9}{10.8}{\rmdefault}{\mddefault}{\updefault}{\color[rgb]{0,0,0}$\disc(x,r',v)$}%
}}}}
\end{picture}%
\caption{Disc $\disc(x,r',v)$.}
\label{chapitre7boulerect}
\end{figure}

We suppose that the event $G(x,r,v)$ happens, and we denote by $U$ a fixed set
satisfying the properties described in the definition of $G(x,r, v)$.
For each $\gamma$ in $\{1/n,...,(\lfloor n\gamma_{\max} \rfloor -1) /n \} $,
we define
$$\left\{ \begin{array}{l} D(\gamma) \,=\, \cyl ( \disc (x, r', v) ,
    \gamma )\,, \\ \p^+ D(\gamma) \,=\, \{y\in D(\gamma) \,|\, \exists z\in
    \ZZ^d_n \,,\,\, (z-x)\cdot v > \gamma \textrm{ and } |z-y|=1 \} \,,\\
\p^- D(\gamma) \,=\, \{y\in D(\gamma) \,|\, \exists z\in
    \ZZ^d_n \,,\,\, (z-x)\cdot v < -\gamma \textrm{ and } |z-y|=1 \}\,.
  \end{array} \right.   $$
These sets are represented in figure \ref{chapitre7boite}.
The sets $\p^+ D(\gamma) \cup\p^- D(\gamma)$ are pairwise disjoint for
different $\gamma$, and we know that
$$ \sum_{\gamma = 1/n,..., (\lfloor n\gamma_{\max} \rfloor -1)/n} \card
((\p^+ D(\gamma) \cap U) \cup (\p^- D(\gamma) \cap U^c)) \,\leq\, 4 \delta
\alpha_{d} r^d n^d   \,,$$
so there exists a $\gamma_0$ in $\{1/n,...,(\lfloor n\gamma_{\max} \rfloor
-1) /n \} $ such that
$$  \card ((\p^+ D(\gamma_0) \cap U) \cup (\p^- D(\gamma_0) \cap U^c))
\,\leq\, \frac{ 4 \delta \alpha_{d} r^d n^d}{\lfloor n\gamma_{\max} \rfloor
  -1} \,\leq\, \frac{5 \delta \alpha_{d} r^d n^{d-1}}{\gamma_{\max}}  $$
for $n$ sufficiently large. We define the event $G^*(x,r,v, \gamma)$
(depending only on the edges in $D(\gamma))$) to be the existence of a set $X
\subset D(\gamma)\cap \ZZ^d_n$ with the following properties:
$$ \left\{ \begin{array}{l} \card ((\p^+ D(\gamma) \cap X) \cup (\p^-
    D(\gamma) \cap X^c))\,\leq\, 5 \delta \alpha_{d} r^d
      n^{d-1}\gamma_{\max}^{-1}\,=\, 5 \delta \alpha_d \rho^{-1} r^{d-1} n^{d-1} \,,\\ V(\p^e X \cap D(\gamma)) \,\leq\,
    \alpha_{d-1} r^{d-1} \nu(v) (1-s) n^{d-1} \,. \end{array}
\right.  $$
We have proved that if $G(x,r, v)$ occurs, there exists a $\gamma$ in $\{
1/n,..., (\lfloor n\gamma_{\max} \rfloor -1)/n  \}$ such that $G^*(x,r, v,
\gamma)$ happens. On $G^*(x,r, v, \gamma)$, we select a set of edges $X$
that satisfies the properties described in the definition of $G^*(B, v(x),
\gamma)$ with a deterministic procedure, and we define
$$\left\{ \begin{array}{l} X^+ \,=\, \{\langle x, y\rangle \,|\, x \in \p^+
    D(\gamma) \cap X \,,\,\, y \notin D(\gamma) \}\,,\\
 X^- \,=\, \{\langle x, y\rangle \,|\, x \in \p^-
    D(\gamma)\smallsetminus X \,,\,\, y \notin D(\gamma) \}\,.
  \end{array} \right. $$
The set of edges $(\p^e X \cap D(\gamma) ) \cup X^+ \cup X^- $ cuts
the top $\partial ^+ D(\gamma_{\max})$ from the bottom $\partial ^-
D(\gamma_{\max})$ of $\C = D(\gamma_{\max})$. If we define
$$ \phi_{\C} \,=\, \phi (\partial ^+ D(\gamma_{\max}) \rightarrow
\partial ^- D(\gamma_{\max}) \textrm{ in } \C)\,, $$
on $G^*(x,r, v, \gamma)$, we have
$$ \phi_{\C} \,\leq\, V(\p^eX \cap D(\gamma)) + V(X^+ \cup X^-) \,.$$
(Recall that $\p^e X \cap D(\gamma)$ is the set of the edges of $\p^e X$
which are included in $D(\gamma)$). Moreover
\begin{align*}
\card( X^+\cup X^-) & \,\leq\, 2d \card((\p^+
    D(\gamma) \cap X)\cup(\p^-
    D(\gamma)\smallsetminus X)) \\ & \,\leq\, 2d  \frac{5 \delta \alpha_{d} r^d
      n^{d-1}}{\gamma_{\max}} \,=\,C  r^{d-1} \delta \rho^{-1} n^{d-1} \,, 
\end{align*}
where $C=10d\alpha_d$ is a constant depending on the dimension. We obtain that
\begin{align*}
\PP [G(x,r,v) ] & \,\leq\, \sum_{\gamma = 1/n ,...,(\lfloor n\gamma_{\max}
  \rfloor -1)/n  } \PP [G^* (x,r,v,\gamma)]\\
& \,\leq\, \sum_{\gamma}
\PP [G^* (x,r,v, \gamma) \cap \{ V(X^+ \cup X^-) \leq \alpha_{d-1}
r^{d-1} \nu(v)n^{d-1} s/4 \}]\\
& \qquad + \PP [G^* (x,r,v, \gamma) \cap \{
V(X^+ \cup X^-) \geq \alpha_{d-1} r^{d-1} \nu(v)n^{d-1} s/4 \}]\,.
\end{align*}
On one hand, we have proved that
\begin{align*}
\PP [G^* (x,r,v, \gamma) & \cap \{ V(X^+ \cup X^-) \leq \alpha_{d-1}
r^{d-1} \nu(v)n^{d-1} s/4 \}] \\ &\,\leq\, \PP  [\phi_{\C} \leq
\alpha_{d-1}r^{d-1}\nu(v) (1-3s/4) n^{d-1}] \,. 
\end{align*}
On the other hand, we have
\begin{align*}
\PP [ &G^* (x,r,v, \gamma)  \cap \{
V(X^+ \cup X^-) \geq \alpha_{d-1} r^{d-1} \nu(v)n^{d-1} s/4 \}]\\
& \,\leq\, \EE \left( \PP( G^* (x,r,v, \gamma) \cap \{ V(X^+ \cup
    X^-) \geq \alpha_{d-1} r^{d-1} \nu(v)n^{d-1} s/4  \} \,|\,
  (t(e))_{e\in D(\gamma)})  \right) \\
& \,\leq\,   \EE \bigg( \PP( G^* (x,r,v, \gamma) \cap \bigcup_{F
   \subset \EE^d_n}( \{X^+\cup X^- =F  \} \\
& \qquad \qquad \cap \{
    V(F) \geq \alpha_{d-1} r^{d-1} \nu(v)n^{d-1} s/4  \} ) \,|\,
  (t(e))_{e\in D(\gamma)})  \bigg)  \\ 
& \,\leq\, \EE \bigg(  \ind{G^* (x,r,v, \gamma)} \sum_{F
    \subset \EE^d_n} \ind{\{X^+\cup X^- =F  \}} \PP (V(F)
  \geq \alpha_{d-1} r^{d-1} \nu(v)n^{d-1} s/4  ) \bigg)\\
&\,\leq\, \PP \left[ \sum_{i=1}^{Cr^{d-1} \delta \rho^{-1} n^{d-1}} t(e_i)
  \geq \alpha_{d-1}  r^{d-1} \nu(v)n^{d-1} s/4   \right]\,, 
\end{align*}
where the last inequality comes from the fact that for all $F$ such that
$\PP [X^+ \cup X^- =F] >0 $, $\card(F) \leq C r^{d-1} \delta \rho^{-1}
n^{d-1}$. Here we have used the following essential
property of $X^+ \cup X^-$: the position of the edges of $X^+ \cup
X^-$ is $\sigma(t(e),e\in D(\gamma))$-measurable, but their capacities are
independent of $(t(e))_{e\in D(\gamma)}$. Finally, we obtain that
\begin{align*}
\PP [G^* (x,r,v, \gamma)] & \,\leq\, \gamma_{\max} n  \PP  [\phi_{\C} \leq
(\alpha_{d-1}r^{d-1}\nu(v)) (1-3s/4) n^{d-1}]\\
& \qquad + \gamma_{\max} n  \PP \left[ \sum_{i=1}^{C r^{d-1}\delta \rho^{-1}
    n^{d-1}} t(e_i) \geq  (\alpha_{d-1}  r^{d-1} \nu(v))n^{d-1} s/4   \right]\,.
\end{align*}
We want to consider cylinders whose basis are hyperrectangles instead of
discs, and the variable $\tau$ instead of $\phi$ in these cylinders,
because we only know the lower large deviations of the flow in this case
(see \cite{RossignolTheret08b}). There exists a constant $c = c(d)$ such
that, for any positive $\kappa$,
there exists a finite family $(A_i)_{i\in I}$ of disjoint closed
hyperrectangles included in $ \disc (x, r', v) $ such that
$$ \left\{ \begin{array}{l} \sum_{i\in I} \H^{d-1}( A_i) \,\geq\,
    \alpha_{d-1} r'^{d-1} - \kappa \,,\\ \sum_{i\in I} \H^{d-2} (\p A_i)
    \leq c r'^{d-2}\,,\end{array}  \right.  $$
(see figure \ref{chapitre7boulerect}).
Thanks to the
max-flow min-cut theorem, we know that for each $i$, the maximal flow
$\tau(\cyl(A_i,\gamma_{max}), v)$ is equal to the smallest capacity of a set
of edges in $\cyl(A_i, \gamma_{\max})$ that cuts the lower half part from
the upper half part of the boundary of the cylinder along the direction
given by $v$. We denote by $\E_i$ such a cutset in $\cyl(A_i,
\gamma_{\max})$. This set of edges is pinned at the boundary of
$A_i$ (which is the common boundary of the two halves of the boundary of
the cylinder $\cyl(A_i, \gamma_{\max})$ between which the flow
$\tau(\cyl(A_i, \gamma_{\max}),v)$ goes). Thus the different sets $\E_i$ in
each cylinder $\cyl(A_i,\gamma_{\max})$ can be glued together along
$\cup_{i\in I} \p A_i$ to create a cutset in $\C$ if we provide some "glue",
i.e., if we add some edges in a
small neighbourhood of $\cup_{i\in I} \p A_i$. For each $i\in I$, we define
the set $\mathcal{P}_i(n) \subset \RR^d$ by
$$  \mathcal{P}_i(n) \,=\, \cyl ( \V(\p A_i, \zeta/n ) \cap \hyp (A_i) ,
\gamma_{\max}) \,, $$
where $\zeta$ is a fixed constant bigger than $2d$, and we denote by
$P_i(n)$ the set of the edges included in $\mathcal{P}_i(n)$. Then
$\cup_{i\in I} E_i \cup P_i(n)$ cuts the top from the bottom of
$\C$. Thanks to the max-flow min-cut theorem again, we thus obtain
that
$$ \phi_{\C} \,\leq\, \sum_{i\in I} \tau(\cyl (A_i, \gamma_{\max}),v) +
V(\cup_{i\in I}P_i(n)) \,. $$
We can evaluate the number of edges in $\cup_{i\in I} P_i(n)$ as follows:
$$ \card(\cup_{i \in I} P_i(n)) \,\leq\, c'
r'^{d-2}\gamma_{\max}n^{d-1} \,\leq\, c' \rho r^{d-1} n^{d-1}   \,,$$
where $c'$ is a constant depending on $\zeta$ and $d$. Therefore
\begin{align*}
 \PP  [\phi_{\C} & \leq \alpha_{d-1}r^{d-1}\nu(v) (1-3s/4) n^{d-1}]\\
& \,\leq\, \PP \left[\sum_{i\in I} \tau(\cyl (A_i, \gamma_{\max}),v ) \leq
\alpha_{d-1}r^{d-1}\nu(v) (1-s/2) n^{d-1}  \right] \\
& \qquad + \PP \left[\sum_{i=1}^{c' \rho r^{d-1} n^{d-1}} t(e_i) \geq
\alpha_{d-1}r^{d-1}\nu(v) \frac{s}{4} n^{d-1}\right]\\
& \,\leq\, \PP \left[\sum_{i\in I} \tau(\cyl (A_i, \gamma_{\max}) ,v) \leq
  (1-s/4) n^{d-1} \sum_{i\in I} \H^{d-1} (A_i ) \nu(v) \right]\\
& \qquad +  \PP \left[\sum_{i=1}^{c' \rho r^{d-1} n^{d-1}} t(e_i) \geq
\alpha_{d-1}r^{d-1}\nu(v) \frac{s}{4} n^{d-1}\right]  \,,
\end{align*}
as soon as the constants satisfy the condition
\begin{equation}
\label{chapitre7cond1}
(\kappa + \alpha_{d-1} (r^{d-1} - r'^{d-1})) (1-s/2) \,\leq\, \sum_{i\in
  I} \H^{d-1} (A_i ) \nu_{\min} s/4 \,.
\end{equation}
Then
\begin{align*}
\PP [G^* (x,r,v, \gamma)] & \,\leq\, \rho r n \sum_{i\in I} \PP  [
\tau(\cyl (A_i, \gamma_{\max}) ,v) \leq \H^{d-1} (A_i ) \nu(v)(1-s/4) n^{d-1} ]\\
& \qquad + \rho r n  \PP \left[ \sum_{i=1}^{C r^{d-1} \delta \rho^{-1}
    n^{d-1}} t(e_i) \geq \alpha_{d-1}  r^{d-1} \nu (v) n^{d-1} s/4   \right]\\
& \qquad + \rho r n \PP \left[\sum_{i=1}^{c' \rho r^{d-1} n^{d-1}} t(e_i) \geq
\alpha_{d-1}r^{d-1}\nu(v)  n^{d-1} s/4 \right]\,. \\
& \,\leq\,  \rho r n \sum_{i\in I} \PP  [
\tau(\cyl (A_i, \gamma_{\max}) ,v) \leq \H^{d-1} (A_i ) \nu(v)(1-s/4) n^{d-1} ]\\
& \qquad + 2 \rho r n  \PP \left[ \sum_{i=1}^{C' (\delta \rho^{-1} + \rho)
    r^{d-1} n^{d-1}}  t(e_i)  \geq   \alpha_{d-1}  r^{d-1}
  \nu(v)n^{d-1} s/2   \right]\,,
\end{align*}
where $C'$ is a constant depending on $\zeta$ and $d$.

%%%%%%%%%%%%%%%%%%%%%%%%%%%%%%%%%%%%%%%%%%%%%%%%%%%%%%%%%%%%%%%%%%%%%%%%%%%%%%%%

\section{Calibration of the constants}

From now on we suppose that the law $\Lambda$ of the capacity of the edges
admits an exponential moment. Then as soon as the constants satisfy the
condition
\begin{equation}
\label{chapitre7cond2}
C' (\rho + \delta \rho^{-1}) r^{d-1} \EE (t(e)) \,<\, (\alpha_{d-1}  r^{d-1}
  \nu_{\min }) \frac{s}{2}\,, 
\end{equation}
the Cram\'er Theorem in $\RR$ allows us to affirm that there exist positive
constants $\D$ and $\D'$ (depending on $\Lambda$, $\delta$, $\rho$, $\zeta$, $s$ and
$d$) such that
$$ \PP \left[ \sum_{i=1}^{C' (\delta \rho^{-1} + \rho)
    r^{d-1} n^{d-1}}  t(e_i)  \geq   (\alpha_{d-1}  r^{d-1}
  \nu(v)n^{d-1} s/2   \right]\, \leq\,\D' e^{-\D n^{d-1}}\,.$$
If we also suppose that $\Lambda(0) <1-p_c(d)$, we know from Theorem
    \ref{thmdevinftau} (Theorem 3.9 in
\cite{RossignolTheret08b}) that there exist a positive
constant $K(d,\Lambda, s)$ and a constant $K'(d, \Lambda,A_i,s)$ such that
$$  \PP  [ \tau(\cyl (A_i, \gamma_{\max}) ,v) \leq \H^{d-1} (A_i )
 \nu(v)(1-s/4) n^{d-1} ] \,\leq\, K' e^{-K n^{d-1} \H^{d-1}(A_i)} \,.$$
We have thus proved that if we can choose, for a fixed $F$, the constants
$\delta$, $\rho$ and $\kappa$ such that for every ball $B$ in the collection of
balls $(B_i)_{i\in I \cup J \cup K}$ the conditions (\ref{chapitre7cond1}) and
(\ref{chapitre7cond2}) are satisfied, then there exists positive constants
$\widetilde{\D}$ and $\hat{\D}$ (depending on $d$, $\Lambda$, $\O$,
$\G^1$, $\G^2$ and $\lambda$) such that
$$ \PP [\phi_n \leq \lambda n^{d-1}] \,\leq\, \hat{\D} e^{-\widetilde{\D} n^{d-1}}
\,,$$
and this yields Theorem \ref{chapitre7devinf}.

We just have to calibrate the constants. In condition (\ref{chapitre7cond2}) appears
the factor $(\rho + \delta \rho^{-1})$: to make it small, we choose $\rho =
\sqrt{\delta}$. Then the condition (\ref{chapitre7cond2}) is equivalent to
$$ \sqrt{\delta} \,<\, \frac{\alpha_{d-1} \nu_{\min} s}{2 C' \EE(t(e))}
\,, $$
for a constant $C'$ that depends on $\zeta$ and $d$, and thus it is
satisfied if we choose $\delta$ small enough (clearly since
$\Lambda(0)<1-p_c(d)$ we know that $\EE(t(e))>0$ and $\nu_{\min}>0$). To see that the condition
(\ref{chapitre7cond1}) can
also be satisfied, we just choose $\kappa \leq \alpha_{d-1}
(r^{d-1}-r'^{d-1})/2$ (so $\kappa$ depends on $\delta$) and we remark that
$$ 1- (\cos \arcsin \sqrt{\delta})^{d-1} \,=\, (d-1) \delta/2 + o(\delta)
\,,$$
so for $\delta$ small enough, condition (\ref{chapitre7cond1}) is satisfied as soon
as
$$ \delta \,\leq\, \frac{2 \nu_{\min} }{12 (d-1) (1-s/2)} \,,$$
which can obviously be satisfied (remember that $s <1$ and $\nu_{\min}>0$). This ends the
proof of Theorem \ref{chapitre7devinf}.

%%%%%%%%%%%%%%%%%%%%%%%%%%%%%%%%%
%\bibliographystyle{plain}
%\bibliography{biblio}

\begin{thebibliography}{10}

\bibitem{ASQU}
P.~Assouad and T.~Quentin~de Gromard.
\newblock Sur la d\'erivation des mesures dans {$\mathbb{R}^n$}.
\newblock 1998.
\newblock Unpublished note.

\bibitem{Bollobas}
B{\'e}la Bollob{\'a}s.
\newblock {\em Graph theory}, volume~63 of {\em Graduate Texts in Mathematics}.
\newblock Springer-Verlag, New York, 1979.
\newblock An introductory course.

\bibitem{Cerf-Pisztora}
Rapha{\"e}l Cerf and {\'A}goston Pisztora.
\newblock Phase coexistence in {I}sing, {P}otts and percolation models.
\newblock {\em Ann. Inst. H. Poincar\'e Probab. Statist.}, 37(6):643--724,
  2001.

\bibitem{CerfTheret09geo}
Rapha\"el Cerf and Marie Th\'eret.
\newblock Law of large numbers for the maximal flow through a domain of
  $\mathbb{R}^d$ in first passage percolation.
\newblock Available from \verb+arxiv.org/abs/0907.5504+.

\bibitem{CerfTheret09sup}
Rapha\"el Cerf and Marie Th\'eret.
\newblock Upper large deviations for the maximal flow through a domain of
  $\mathbb{R}^d$ in first passage percolation.
\newblock Available from \verb+arxiv.org/abs/0907.5499+.

\bibitem{Cerf:StFlour}
Raphaël Cerf.
\newblock The {W}ulff crystal in {I}sing and percolation models.
\newblock In {\em \'Ecole d'\'Et\'e de Probabilit\'es de Saint Flour}, number
  1878 in Lecture Notes in Mathematics. Springer-Verlag, 2006.

\bibitem{FAL}
K.~J. Falconer.
\newblock {\em The geometry of fractal sets}, volume~85 of {\em Cambridge
  Tracts in Mathematics}.
\newblock Cambridge University Press, Cambridge, 1986.

\bibitem{Kesten:StFlour}
Harry Kesten.
\newblock Aspects of first passage percolation.
\newblock In {\em \'Ecole d'\'Et\'e de Probabilit\'es de Saint Flour XIV},
  number 1180 in Lecture Notes in Mathematics. Springer-Verlag, 1984.

\bibitem{Kesten:flows}
Harry Kesten.
\newblock Surfaces with minimal random weights and maximal flows: a higher
  dimensional version of first-passage percolation.
\newblock {\em Illinois Journal of Mathematics}, 31(1):99--166, 1987.

\bibitem{RossignolTheret08b}
R.~Rossignol and M.~Th\'eret.
\newblock Lower large deviations and laws of large numbers for maximal flows
  through a box in first passage percolation.
\newblock Available from \verb+arxiv.org/abs/0801.0967v2+, 2009.

\bibitem{Zhang}
Yu~Zhang.
\newblock Critical behavior for maximal flows on the cubic lattice.
\newblock {\em Journal of Statistical Physics}, 98(3-4):799--811, 2000.

\bibitem{Zhang07}
Yu~Zhang.
\newblock Limit theorems for maximum flows on a lattice.
\newblock Available from \verb+arxiv.org/abs/0710.4589+, 2007.

\end{thebibliography}

\def\cprime{$'$}

\end{document}